\documentclass[a4paper,11pt]{article}

\usepackage{amsthm}
\usepackage{hyperref}
\usepackage{amssymb}
\usepackage{latexsym}
\usepackage{amsmath}
\usepackage{amsfonts}
\usepackage{version}
\usepackage[dvips]{graphicx}
\usepackage[latin1]{inputenc}

\usepackage[LGR,T1]{fontenc}%
\usepackage[francais,english]{babel}
\usepackage{url}
\usepackage{enumerate}
\usepackage{hyperref}
\usepackage{color}

\usepackage{fancyhdr}
\newtheorem{Theorem}{Theorem}[part]
\newtheorem{Definition}{Definition}[part]
\newtheorem{Proposition}{Proposition}[part]

\newtheorem{Corollary}{Corollary}[part]
\newtheorem{Remark}{Remark}[part]
\newtheorem{Example}{Example}[part]

\makeatletter \@addtoreset{equation}{section}

\@addtoreset{Definition}{section}

\@addtoreset{Theorem}{section}

\@addtoreset{Proposition}{section}

\@addtoreset{Property}{section}

\@addtoreset{Assumption}{section}

\@addtoreset{Corollary}{section}

\@addtoreset{Lemma}{section}

\@addtoreset{Remark}{section}

\@addtoreset{Example}{section}

\def \Prod{\displaystyle\prod}

\def\={\;=\;}

\def\be{\begin{eqnarray}}
\def\ee{\end{eqnarray}}
\def\b*{\begin{eqnarray*}}
\def\e*{\end{eqnarray*}}
\def\beq{\begin{equation}}
\def\eeq{\end{equation}}

\def\1{{\bf 1}}

\def \proof{{\noindent \bf Proof. }}
\def \eproof{\hbox{ }\hfill$\Box$}

\def \R{\mathbb{R}}

\def\P{\mathbb{P}}

\def\cA{{\cal A}}
\def\cB{{\cal B}}

\def\cE{{\cal E}}
\def\cF{{\cal F}}
\def\cG{{\cal G}}
\def\cH{{\cal H}}

\def\cM{{\cal M}}

\def\cS{{\cal S}}
\def\cT{{\cal T}}

\def\ti{{t_i}}

\newcommand{\ud}
    {\mathrm{d}}

\newcommand{\esp}[1]
    {\ensuremath{%
     \mathbb{E}\!\!\left[#1\right]}}

\newcommand{\EFp}[2]
    {\ensuremath{
     \mathbb{E}_{#1}\!\!\left[#2\right] }}

\newcommand{\HP}[1] 
    {\ensuremath{\mathcal{H}^{#1}}}

\newcommand{\HPb}[1] 
    {\ensuremath{\mathcal{H}_\beta^{#1}}}

\newcommand{\Eb} 
    {\ensuremath{E_\beta}}

\newcommand{\ESP}[1] 
    {\ensuremath{\mathbb{S}^{#1}}}

\newcommand{\SPb}[1] 
    {\ensuremath{\mathcal{S}_\beta^{#1}}}

\newcommand{\Pot}[1] 
    {\ensuremath{\Pi_c^{#1}}}

\newcommand{\DM}[1] 
    {\ensuremath{\mathbb{D}^{1,#1}}}

\newcommand{\LM}[1] 
    {\ensuremath{\mathbb{L}_a^{1,#1}}}

\newcommand{\Lp}[1] 
    {\ensuremath{ \mathbf{L}^{#1} }}

\newcommand{\NH}[2] 
    { \ensuremath{ ||#2||_{\mathcal{H}^{#1}} } }

\newcommand{\NHb}[2] 
    { \ensuremath{ ||#2||_{\mathcal{H}_\beta^{#1}} } }

\newcommand{\set}[1]
    {\ensuremath{\{ #1 \}}}

\newcommand{\HYP}[1]
    {\ensuremath{({\bold H#1} ) }}

\newcommand{\Zr}{\ensuremath{ Z^{d\Re} }}

\newcommand{\ZriO}[1]{\ensuremath{ (Z^{d\Re})^{.i} }}

\newcommand{\ZrkiO}[1]{\ensuremath{ (Z^{d\Re})^{ki} }}

\newcommand{\tYr}{\ensuremath{\widetilde Y^{d\Re} }}

\newcommand{\tZr}{\ensuremath{\widetilde Z^{d\Re} }}

\newcommand{\Xp} 
    {\ensuremath{X^{\pi} }}
\newcommand{\Xpt}[1] 
    {\ensuremath{X^{\pi}_{t_{#1}} }}

\newcommand{\Xpr}[1] 
    {\ensuremath{X^{\pi}_{r_{#1}} }}

\newcommand{\Ypb} 
    {\ensuremath{Y^{\pi,b} }}

\newcommand{\Ylb} 
    {\ensuremath{\tYr }}
\newcommand{\Ylbi} 
    {\ensuremath{\tilde Y^{b,i} }}
\newcommand{\Zlb} 
    {\ensuremath{\tZr }}

\newcommand{\dY} 
    {\ensuremath{\delta Y }}

\newcommand{\Ypt}[1] 
    {\ensuremath{Y^{\pi}_{t_{#1}} }}

\newcommand{\Ypr}[1] 
    {\ensuremath{Y^{\pi}_{r_{#1}} }}

\newcommand{\Ylp} 
    {\ensuremath{\widetilde{Y}^{\pi} }}

\newcommand{\Ylpt}[1] 
    {\ensuremath{\widetilde{Y}^{\pi}_{t_{#1}} }}

\newcommand{\dZ}
    {\ensuremath{\delta Z }}

\newcommand{\Zlp} 
    {\ensuremath{\widetilde{Z}^{\pi} }}

\newcommand{\Zpt}[1] 
    {\ensuremath{Z^{\pi}_{t_{#1}} }}

\newcommand{\Zpr}[1] 
    {\ensuremath{Z^{\pi}_{r_{#1}} }}

\newcommand{\Zlpt}[1] 
    {\ensuremath{\widetilde{Z}^{\pi}_{t_{#1}} }}

\newcommand{\Zhpt}[1] 
    {\ensuremath{\widehat{Z}^{\pi}_{t_{#1}} }}

\newcommand{\Zbb} 
    {\ensuremath{\bar \Zr }}

\newcommand{\wh}[1]{\ensuremath{ \widehat{#1} }}

\newcommand{\Yi}{\ensuremath{ Y^{i} }}

\newcommand{\Zi}{\ensuremath{ Z^{i} }}

\renewcommand{\Xi}[1]{X_{i #1}}
\renewcommand{\ti}[1]{t_{i #1}}
\newcommand{\tk}[1]{t_{k #1}}

\renewcommand{\Yi}[1]{Y_{i #1}}
\renewcommand{\Zi}[1]{Z_{i #1}}

\newcommand{\whyi}[1]{\widehat{y}_{i #1}}
\newcommand{\whzi}[1]{\widehat{z}_{i #1}}

\newcommand{\tYi}[1]{\tilde{Y}_{i #1}}
\newcommand{\tZi}[1]{\tilde{Z}_{i #1}}
\newcommand{\wtYi}[1]{\widetilde{Y}_{i #1}}
\newcommand{\wtZi}[1]{\widetilde{Z}_{i #1}}
\newcommand{\whYi}[1]{\widehat{Y}_{i #1}}
\newcommand{\whZi}[1]{\widehat{Z}_{i #1}}
\newcommand{\dYi}[1]{\delta Y_{i #1}}
\newcommand{\dZi}[1]{\delta Z_{i #1}}
\newcommand{\dYj}[1]{\delta Y_{j #1}}

\newcommand{\dYk}[1]{\delta Y_{k #1}}
\newcommand{\dZk}[1]{\delta Z_{k #1}}

\newcommand{\Ui}[1]{U_{i #1}}

\newcommand{\dUj}[1]{\delta U_{j #1}}

\newcommand{\hi}{h_{i}}
\newcommand{\errYi}[1]{\zeta^{Y}_{i #1}}
\newcommand{\errZi}[1]{\zeta^{Z}_{i #1}}

\newcommand{\errYk}[1]{\zeta^{Y}_{k #1}}
\newcommand{\errZk}[1]{\zeta^{Z}_{k #1}}

\newcommand{\phiy}{\varphi^Y}
\newcommand{\phiz}{\varphi^Z}
\newcommand{\Wxi}[1]{W^{x,\ti{}}_{\ti{#1}}}
\newcommand{\Woi}[1]{W^{0,\ti{}}_{\ti{#1}}}
\newcommand{\whWxi}[1]{\widehat{W}^{x,\ti{}}_{\ti{#1}}}
\newcommand{\whWoi}[1]{\widehat{W}^{0,\ti{}}_{\ti{#1}}}
\newcommand{\Wpi}[1]{\widehat{W}^{x,0}_{\ti{#1}}} 
\newcommand{\Ryi}[1]{R^{Y}_{i,#1}}
\newcommand{\Rzi}[1]{R^{Z}_{i,#1}}
\newcommand{\whRyi}[1]{\widehat{R}^{Y}_{i,#1}}
\newcommand{\whRzi}[1]{\widehat{R}^{Z}_{i,#1}}
\newcommand{\whHi}[1]{\widehat{H}_{i#1}}
\newcommand{\whEFp}[2]
    {\ensuremath{%
     \widehat{\mathbb{E}}_{#1}\!\left[#2\right] }}
\newcommand{\whesp}[1]
    {\ensuremath{%
     \widehat{\mathbb{E}}\!\left[#1\right] }}

\title{\vspace{-5mm} Linear multi-step schemes for BSDEs}

\author{Jean-Fran\c{c}ois CHASSAGNEUX\footnote{Departement of Mathematics, Imperial College London.  { \sf j.chassagneux@imperial.ac.uk} }
              }
\date{This version: October 2012. (submitted)}
\begin{document}
\maketitle

\begin{abstract}
We study the convergence rate of a class  of linear multi-step methods for BSDEs.  We show that, under a sufficient condition on the coefficients, the schemes enjoy a fundamental stability property. Coupling this result to an analysis of the truncation error allows us to design approximation with arbitrary order of convergence. Contrary to the analysis performed in \cite{zhazha10}, we consider general diffusion model and  BSDEs with driver depending on $z$.
The class of methods we consider contains well known methods from the ODE framework as Nystrom, Milne or Adams methods. We also study a class of Predictor-Correctot methods based on Adams methods. Finally, we provide a numerical illustration of the convergence of some methods.
\end{abstract}

\vspace{2mm}

\vspace{1mm}

\noindent {\bf MSC Classification (2000):} 60H10, 65C30.



\section{Introduction}

In this paper, we are interested in the discrete-time approximation of solutions of (decoupled) Backward Stochastic Differential Equation (BSDE), i.e$.$ a triplet $(X,Y,Z)$ satisfying
%
\begin{align}
X_{t}& = X_{0} + \int_{0}^{t}b(X_{s})\ud s + \int_{0}^{t}\sigma(X_{s})\ud W_{s},\label{eq sde X}
\\
Y_{t} &= g(X_{T})+ \int_{t}^{T} f(Y_{t},Z_{t}) \ud t  -\int_{t}^{T} Z_{t} \ud W_{t} \,.
\label{eq bsde YZ}
\end{align}
The function $(b,\sigma): \R^d \mapsto \R^d \times \cM_d$,
and $f:\R \times \R^d \mapsto \R$ are Lipschitz-continuous function,  $g: \R^d \mapsto \R$ is differentiable with continuous and bounded first derivative\footnote{These assumptions will be strengthened in the following sections.}.   The positive constant $T$ is given and $W$ is a Brownian motion supported by a filtered probability space $(\Omega, \cF, (\cF)_{0 \le t \le T}, \P)$. The process $Y$ is a one-dimensional stochastic
process, the processes $X$  and $Z$ are valued in $\R^d$ and $Z$ is  written, 
by convention, as a row vector. Under the Lipschitz assumption on the coefficients, the processes $X$ and $Y$ belong to the set $\cS^2$ of continuous adapted processes with square integrable supremum and $Z$ belongs to $\cH^2$, the set of progressively measurable processes satisfying $\esp{\int_0^T |Z_s|^2 \ud s}$.

\vspace{1mm}
The existence and uniqueness of solutions of  the system \eqref{eq sde X} -\eqref{eq bsde YZ} was first
addressed by Pardoux and Peng in \cite{parpen90}. 
Moreover,  in \cite{parpen92}, they show that
\[
Y_t=u(t,X_t),\ \ \  \ Z_t=\nabla u^{\!\top}\!(t,X_t)\sigma(X_t),  \ \ \ \ t\in [0,T],
\]
where $u\in C^{1,2}([0,T]\times \mathbb{R}^d)$ is the solution of the final 
value Cauchy problem 
\begin{align}
L^{(0)}u(t,x)&=-f\left( u(t,x),\nabla u^{\!\top}\!(t,x)\sigma\left( x\right) \right)
,\quad t\in \lbrack 0,T),\,x\in \mathbb{R}^{d} \label{b_pde_t}
\\
u(T,x)&=g (x),\quad x\in \mathbb{R}^{d} \label{b_pde_T}
\end{align}
with $L^{(0)}$  defined to be the second order differential operator\vspace{-3mm} 
\begin{equation}
L^{(0)}=\partial _{t} + \sum_{i=1}^{d}b_i\partial_{x_i}+ \frac12 \sum_{i,j=1}^{d}a_{ij}
\partial_{x_i}\partial_{x_j}\label{operator},
\end{equation}
and $a=a_{ij}=\sigma\sigma^{\top}$.

\vspace{1mm}
To approximate \eqref{eq sde X}-\eqref{eq bsde YZ}, one has to come up with an approximation of the SDE part  and the BSDE part.
Obtaining approximations of the distribution of the forward component $X$ has been largely resolved in the last thirty years. There is a large literature on the subject and one can refer to  
\cite{klopla92}  and the references therein for a systematic study of numerical methods  for approximating $X$. 

Here, we  focus on the approximation of $(Y,Z)$ instead.
Numerical methods  approximating this backward component  have already been proposed. They are mainly based on a Euler approximation, see \cite{boutou04, zha04, goblab07, criman09} and the references therein. These methods have been successfully extended to a broader class of BSDEs: reflected BSDEs \cite{balpag03, cha08}, BSDEs with jumps \cite{boueli08}, BSDEs with driver of quadratic growth \cite{ric10}, see also the reference therein. In a very specific framework, \cite{zhaliy12, zhache06, zhawan09, zhazha10} proposed  some high order methods to approximate the solution of the BSDE. Recently high order method of Runge-Kutta type have been studied \cite{criman10, chacri12} in the general framework of  \eqref{eq sde X}-\eqref{eq bsde YZ}.

In this paper, we consider another type of high order method, very well known for ODEs, namely linear multi-step  methods.  

\vspace{1mm}
The approximations presented below are associated to an arbitrary, but fixed, 
partition $\pi$ of the interval $[0,T]$, $\pi = \set{t_{0} = 0 <\dots<t_i<t_{i+1}<\dots< t_{n} = T}$. We define $\hi = \ti{+1}-\ti{}$, $i=0,...,n-1$ and $|\pi| = \max_{i}\hi$ and denote by $(\Yi{},\Zi{})$ the approximation of $(Y_{t_i},Z_{t_i})$ for $i=1,...,n$. The construction of the approximating process is done in a recursive manner, backwards in time. We describe in the following the salient features of the class of approximations considered in this paper.

\begin{Definition} (Linear multi-step methods)\label{de multi-step scheme}

(i) To initialise the scheme with $r$ steps, $r\ge 1$, we are given  $r$ terminal condition $(Y_{n-j},Z_{n-j})$,  $\cF_{t_{n-j}}$-measurable square integrable random variables, $0 \le j \le r-1$.

(ii) For $i \le n-r$, the  computation of $(Y_i,Z_i)$ involves $r$ steps  and is given by
\begin{align*} 
\left \{
\begin{array}{rcl}
\Yi{} &=& \EFp{\ti{}}{ \sum_{j=1}^r a_{j}\Yi{+j} + h \sum_{j=0}^r b_{i,j} f(\Yi{+j},\Zi{+j}) }
\vspace{2mm}
\\
\Zi{} &=&\EFp{\ti{}}{ 
		\sum_{j = 1}^r \alpha_{j} H^Y_{i,j} \Yi{+j} + h \sum_{j=1}^r \beta_{i,j} H^f_{i,j} f(\Yi{+j},\Zi{+j}) 
}
\end{array}
\right.
\end{align*}

where $a_{j}$, $b_{i,j}$, $\alpha_{j}$, $\beta_{i,j}$ are real numbers satisfying
\begin{align*}
|a_{j}|+ |b_{i,j}| + |\alpha_{j}|+ |\beta_{i,j}|  \le \Lambda \;, \; 0 \le i \le n-r\;, 0 \le j \le r\;,
\end{align*}
and  $\Lambda$ is a positive constant.
We impose the so-called pre-consistency condition  i.e.
\begin{align*}
\sum_{j=1}^r{a_j} = \sum_{j=1}^r{\alpha_j} = 1\;.
\end{align*}

The coefficients $H^Y_{i,j}$, $H^f_{i,j}$, $0 \le i \le n-r$, $1 \le j \le r$ are $\cF_{\ti{+j}}$-measurable random variables satisfying, for all $j$,
\begin{align*}
h_i\esp{|H^Y_{i,j}|^2+|H^f_{i,j}|^2} \le \Lambda \quad\text{ and } \quad \EFp{\ti{}}{H^Y_{i,j}} = \EFp{\ti{}}{H^f_{i,j}} = 0 \,.
\end{align*}

\end{Definition}

\vspace{2mm}
\begin{Remark} \label{re intro}
(i) The value $(Y_n,Z_n)$ is generally given by $(g(X_T),\nabla g^\top (X_T) \sigma(X_T))$. If $r>1$, one needs to specify other initialisation values. This choice is important because it will impact the global rate of convergence. One can use Runge-Kutta type scheme \cite{chacri12} with high order of convergence.

(ii) When $r=1$, schemes \ref{de multi-step scheme}  are one-step scheme. See \cite{boutou04, zha04, goblab07, criman10, chacri12} and the references therein for a study of these schemes.
\end{Remark}

\vspace{2mm}
The global error we investigate here is a \emph{time discretization error} and is, given a grid $\pi$,  $(\cE_Y(\pi),\cE_Z(\pi))$ with
\begin{align*}
\cE_Y(\pi) := \max_i \esp{|Y_{\ti{}}-Y_i|^2} \text{ and }  \cE_Z(\pi) := \sum_i \hi \esp{|Z_{\ti{}}-Z_i|^2}.
\end{align*}
%

To implement high order scheme in practice, we need to specify a particular form for the $H$-coefficient appearing in Definition \ref{de multi-step scheme} above. Let us first introduce a special class of random variables, which was already considered in \cite{chacri12}.

\begin{Definition}\label{de psi-H} (i) For $m \ge 0$, we denote by $\cB^{m}_{[0,1]}$ the set of bounded measurable  function
 $\psi :[0,1] \rightarrow \R$   satisfying
\begin{align*}
\int_{0}^{1}\psi(u)\ud u = 1 \text{ and if } m \ge 1,\;  \int_{0}^{1}\psi(u)u^{k}\ud u =0\;,\; 1\le k \le m.
\end{align*} 

(ii) Let $(\psi^{\ell})_{1\le \ell \le d} \in \cB^{m}_{[0,1]}$,   
for $t \in [0,T]$ and $h>0$ s.t. $t+h\le T$, we define,
\begin{align*}
H^{\psi}_{t,h} := (\frac1h\int_{t}^{t+h}\psi^{\ell}(\frac{u-t}{h}) \ud W^\ell_{u})_{1\le \ell \le d} \;,
\end{align*}
which is a row vector.

By convention, we set $H^{\psi}_{t,0} = 0$.
\end{Definition}

\vspace{2mm}
In the sequel, when studying the order of convergence of the scheme and depending of the order we want to retrieve, we will  assume that, for $1\le j \le r$, 
\begin{align}
H^Y_{i,j} := H^\psi_{\ti{},jh} \text{ and } H^f_{i,j} := H^\phi_{\ti{},jh} \label{eq def H scheme}
\end{align}
for some functions $\psi$ and $\phi$ in $\cB^m$, $m \ge 0$, see Theorem \ref{th main conv res}  below.

\vspace{2mm}
The convergence analysis is done in a classical way. We first prove a fundamental stability property for the schemes, under a reasonable sufficient condition, see Proposition \ref{pr L2 stab}. Then, assuming  smoothness of the value function $u$ given by \eqref{b_pde_t}-\eqref{b_pde_T}, we study the truncation error associated to the above methods. We prove a sufficient condition on the coefficient to retrieve methods of any order. These two steps allow us to retrieve general convergence and design new high order method for BSDEs.  Contrary to the analysis performed in \cite{zhazha10}, we work with general diffusion model given by \eqref{eq sde X} and  BSDEs with driver depending on $z$. 
As an example of application, we extend some classical scheme used in the ODE framework and then proceed with the study of  Adams type methods. Based on these methods, we also design Predictor-Corrector methods and study their convergence. To the best of our knowledge, it is the first time that these methods are considered for BSDEs. 
Finally, we illustrate our theoretical results with some numerical experiments showing empirical convergence rates. 

\vspace{2mm}
The rest of this paper is organised as follows. In section 2, we prove our  general convergence result which relies heavily on a stability property. In section 3, we study Adams methods  and Predictor-Corrector methods in the context of BSDEs. The main results are stated in the multi-dimensional case but for the reader's convenience the proofs are done with $d=1$. Finally, in section 4, we provide a numerical example.

\vspace{4mm}
\paragraph{Notations}
We denote by $\mathcal{M}_d$ the set of matrices with $d$ lines and $d$ columns.
For a matrix $A\in \cM_d$, $\text{Tr}[A]$ denotes its trace, $A^{.j}$ its $j$-th column, $A^{i.}$ its $i$-th row,  and $A^{ij}$ the $i$-th term of  $A^{.j}$. $I_d$ is the identity matrix of $\cM_d$.
The transpose of a matrix or a vector $y$ will be denoted $y^\top$. 
 The sup-norm  for both vectors and matrix is denoted $|\,.\,|_\infty$.
 
 In the sequel $C$ is a positive constant whose value may change from line to line depending on $T$, $d$, $\Lambda$, $X_0$ but which does not depend on $\pi$. We write $C_p$ if it depends on some positive parameters $p$.

For $t \in \pi$, $R$ a random variable and $r$ a real number, the notation $R=O_t(r)$ 
means that $|U| \le \lambda^\pi_t u$ where $\lambda^\pi_t$ is a positive random variable satisfying:
\begin{align*}
\esp{|\lambda^\pi_t|^p} \le C_p\;,
\end{align*}
for all $p > 0$, $t \in \pi$ and all grid $\pi$.

%
%
%
%
%


\section{General convergence results}

In this part, we study the convergence properties of the schemes given in Definition \ref{de multi-step scheme}.

We first establish a stability property for the schemes.   We then state a sufficient condition on the coefficients which allows us to retrieve high order schemes.

\subsection{$L^2$-stability}

%
%
%

To investigate the \emph{stability} of the schemes given in Definition \ref{de multi-step scheme}, we introduce a pertubed scheme
\begin{align} \label{eq multi-step scheme stab pert}
\left \{
\begin{array}{rcl}
\tYi{} &=& \EFp{\ti{}}{ \sum_{j=1}^r a_{j}\tYi{+j} + h \sum_{j=0}^r b_{i,j} f(\tYi{+j},\tZi{+j})} +  \errYi{} 
\vspace{1mm}
\\
\tZi{} &=&\EFp{\ti{}}{ 
		\sum_{j = 1}^r \alpha_{j} H^Y_{i,j} \tYi{+j} + h \sum_{j=1}^r \beta_{i,j} H^f_{i,j} f(\tYi{+j},\tZi{+j})  } + \errZi{}
\end{array}
\right.
\end{align}

where $ \errYi{}$, $\errZi{}$ are random variables belonging to $L^2(\cF_{\ti{}})$, for $i \le n-r$.

\vspace{2mm}
The notion of \emph{stablity} we consider here is the following.

\begin{Definition}($L^2$-Stability)\label{de stab}
The scheme given in Definition \ref{de multi-step scheme} is said to be $L^2$-stable if 
\begin{align*}
\max_{0 \le i \le n-r} \esp{|\dYi{}|^2} & + \sum_{i=0}^{n-r} \hi \esp{|\dZi{}|^2}
\le  \\
& C \Big(\max_{0\le j \le r-1} \esp{|\dY_{n-j}|^2 + |\pi| |\dZ_{n-j}|^2} + |\pi|\sum_{i=0}^{n-r} \esp{\frac{1}{\hi^2} |{\errYi{}}|^2 +  |{\errZi{}}|^2} \Big)
\end{align*}

for all sequences $\errYi{}$,$\errZi{}$ of $L^2(\cF_{\ti{}})$-random variable, $i\le n-r$, and terminal values $(Y_{n-j}, Z_{n-j})$, $(\tilde{Y}_{n-j}, \tilde{Z}_{n-j})$ belonging to $L^2(\cF_{t_{n-j}})$, $0\le j \le r-1$. 
\end{Definition}


%
\vspace{2mm}

\begin{Proposition} 
\label{pr L2 stab}
Assume that the following holds
\begin{itemize}
\item[\HYP{c}]  The coefficients $(a_j)$ are non-negative, $\sum_{j=1}^r a_j = 1$ and for $1 \le j \le r$, $a_j=0\implies \alpha_j = 0$,
\end{itemize}
then, the scheme given in Definition \ref{de multi-step scheme}  is $L^2$-stable, recalling Definition \ref{de stab}. 
\end{Proposition}

\proof

We define $U_i = (Y_i,\dots,Y_{i+r-1})^\top$, $\tilde{U}_i = (\tYi{},\dots,\tYi{+r-1})^\top$, 
and

 \begin{align*}
 \Phi^Y_i &=  \left( \begin{array}{c}
\sum_{j=0}^r b_{i,j} f(\Yi{+j},\Zi{+j}) \\
\mathbf{0} \\
 \end{array} \right)_{r,1} 
 \text{ , } 
  \tilde{\Phi}^Y_i =  \left( \begin{array}{c}
\sum_{j=0}^r b_{i,j} f(\tYi{+j},\tZi{+j}) \\
\mathbf{0} \\
 \end{array} \right)_{r,1} 
 \\
  \text{ and } &
  \tilde{\Theta}^Y_i =  \left( \begin{array}{c}
\errYi{}\\
\mathbf{0} \\
 \end{array} \right)_{r,1} 
 \end{align*}
 and
 denote  $\delta U_i = U_i-\tilde{U}_i$, $\delta \Phi^Y_i = \Phi^Y_i - \tilde{\Phi}^Y_i$,
 \begin{align*}
a =(a_1, \dots,a_r),\; \alpha=(\alpha_1, \dots,\alpha_r) \text{ and } A =  \left( \begin{array}{c|c}
a_1, \dots, a_{r-1} & a_r \\
 \hline
 I_{r-1} & \mathbf{0}  \\
 \end{array} \right)_{r,r} \;.
\end{align*}

The scheme and the pertubed scheme rewrite then for the $Y$ part
\begin{align*}
\EFp{\ti{}}{\Ui{}} &= \EFp{\ti{}}{ A \Ui{+1} + h \Phi^Y_i }
\end{align*}
and
\begin{align*}
\EFp{\ti{}}{\tilde{U}_i} &= \EFp{\ti{}}{ A\tilde{U}_{i+1} + h \tilde{\Phi}^Y_i + \Theta^Y_i } \;.
\end{align*}

1.a  

For $i \le j \le n-r$, we  compute that
\begin{align*}
|\EFp{\ti{}}{\dUj{}}|_\infty \le | A |_\infty |\EFp{\ti{}}{\dUj{+1}}|_\infty + h_j| \EFp{\ti{}}{ \delta \Phi^Y_j}|_\infty + |\EFp{\ti{}}{\Theta^Y_j}|_\infty
\end{align*}
Under $\HYP{c}$, we observe that $|A|_\infty = 1$ and we get
\begin{align*}
| \EFp{\ti{}}{\dUj{}}|_\infty \le  |\EFp{\ti{}}{\dUj{+1}}|_\infty + h_j|\EFp{\ti{}}{ \delta \Phi^Y_j}|_\infty + |\EFp{\ti{}}{\Theta^Y_j}|_\infty
\end{align*}
Iterating on $j$, we compute that
\begin{align*}
| \EFp{\ti{}}{\dUj{}}|_\infty \le  |\EFp{\ti{}}{\delta U_{n-r+1}}|_\infty + \sum_{k=j}^{n-r}h_k|\EFp{\ti{}}{ \delta \Phi^Y_k}|_\infty + \sum_{k=j}^{n-r} |\EFp{\ti{}}{\Theta^Y_k}|_\infty \;.
\end{align*}


In particular, we have for $i = j$, and $|\pi|$ small enough,
\begin{align}\label{eq stab Y step 1}
|\dYi{}| \le C\Big( \sum_{k=i}^{n-r}h_k\EFp{\ti{}}{|\dYk{}| + |\dZk{}|}+ \sum_{k=i}^{n-r}\EFp{\ti{}}{|\errYk{}|} 
		+ \sum_{k=n-r+1}^{n}\EFp{\ti{}}{|\dY_{k}|} 
		 \Big)\;.
\end{align}

We then compute
\begin{align}\label{eq stab Y step 2}
\esp{|\dYi{}|^2} \le C\Big( |\pi|\sum_{k=i}^{n-r}\esp{|\dYk{}|^2} + \sum_{k=i}^{n-r}h_k\esp{|\dZk{}|^2} &+ 
\sum_{k=i}^{n-r}\frac1{h_k}\esp{|\errYk{}|^2} 
		+ \sum_{k=n-r+1}^{n}\esp{|\dY_{k}|^2} \Big) \;.
\end{align}

1.b We will now control the term $ h\sum_{k=i}^{n-r}\esp{|\dZk{}|^2}$ appearing in \eqref{eq stab Y step 2}. 

Using Cauchy-Schwartz inequality, we obtain  that, if $a_j \neq 0$ then
\begin{align*}
|\EFp{\ti{}}{\alpha_j H_{i,j}^Y\dYi{+j}}|^2 
& \le C(a_j \EFp{\ti{}}{|\dYi{+j}|^2 - a_j |\EFp{\ti{}}{\dYi{+j}}|^2})
\end{align*}
which leads to, under \HYP{c},
\begin{align} \label{eq stab Z step 1}
h_i\esp{|\dZi{}|^2} &\le C\Big(\sum_{j=1}^r a_j\esp{|\dYi{+j}|^2} - \sum_{j=1}^r\esp{a_j|\EFp{\ti{}}{ \dYi{+j}}|^2})
\nonumber
\\
&\quad+ |\pi|^2 \sum_{j=1}^r\esp{|\dYi{+j}|^2 + |\dZi{+j}|^2} + |\pi|\esp{|\errZi{}|^2} \;.
\Big) 
\end{align}
Under $\HYP{c}$, we have that, 
\begin{align*}
- \sum_{j=1}^r\esp{a_j |\EFp{\ti{}}{ \dYi{+j}}|^2} & \le -\esp{|\sum_{j=1}^r\EFp{\ti{}}{a_j \dYi{+j}}|^2} \;.
\end{align*}

Then, recalling that
\begin{align*}
\sum_{j=1}^r\EFp{\ti{}}{a_j \dYi{+j}} = \delta Y_i - h_i\sum_{j=0}^r \EFp{\ti{}}{\delta \Phi^Y_{j+r}} - \errYi{}
\end{align*}
we compute 
\begin{align*}
- \sum_{j=1}^r\esp{|\EFp{\ti{}}{a_j \dYi{+j}}|^2}  &\le -\esp{| \dYi{}|^2} + C|\pi|\esp{| \dYi{}|\sum_{j=0}^r\EFp{\ti{}}{|\dYi{+j}|+|\dZi{+j}|}} + C\esp{| \dYi{}|\,|\errYi{}| }
\end{align*}
which leads, for $0<\epsilon \le 1$ to be fixed later on, to
\begin{align*}
- \sum_{j=1}^r\esp{|\EFp{\ti{}}{a_j \dYi{+j}}|^2}  &\le -\esp{|\dYi{}|^2} + C |\pi| \sum_{j=0}^r\EFp{\ti{}}{\frac1\epsilon|\dYi{+j}|^2+\epsilon|\dZi{+j}|^2} +  \frac{C}{h_i}\esp{|\errYi{}|^2}\;.
\end{align*}

Combining the last inequality with \eqref{eq stab Z step 1} and summing over $i$, we obtain, for $|\pi|$ small enough
\begin{align*}
\sum_{k=i}^{n-r} h_k\esp{|\dZk{}|^2} &\le C\Big(\sum_{k=i}^{n-r}(\sum_{j=1}^r a_j\esp{|\dYi{+j}|^2} - \esp{|\dYi{}|^2}) + C(1+\frac1\epsilon) |\pi| \sum_{k=n-r+1}^{n}\esp{|\dYk{}|^2+|\dZk{}|^2} 
\\
&+ C(1+\frac1\epsilon)|\pi| \sum_{k=i}^{n-r}\esp{|\dYk{}|^2} +\frac{C}{\epsilon} |\pi| \sum_{k=i}^{n-r}\esp{|\dZk{}|^2} +|\pi|\sum_{k=i}^{n-r}\esp{|{\errZk{}}|^2}+ \sum_{k=i}^{n-r} \frac{C}{h_k}\esp{|{\errYk{}}|^2}
\Big) 
\end{align*}

Using \HYP{c}, setting $\epsilon := \frac{C}{2}$, we then obtain
\begin{align}
\sum_{k=i}^{n-r} h_k \esp{|\dZk{}|^2} &\le C\Big( \sum_{k=n-r+1}^{n}\esp{|\dYk{}|^2+|\pi||\dZk{}|^2} 
+  |\pi| \sum_{k=i}^{n-r}\esp{|\dYk{}|^2} + \sum_{k=i}^{n-r}\esp{\frac1{h_k}|{\errYk{}}|^2 +h_k|{\errZk{}}|^2}
\Big) \label{eq stab interm Z}
\end{align}

1.c Combining the last inequality with \eqref{eq stab Y step 2}, we get
\begin{align}\label{eq stab Y step 1.c-1}
\esp{|\dYi{}|^2} \le C\Big( |\pi|\sum_{j=i}^{n-r}\esp{|\dYj{}|^2} 
+ 
\sum_{k=i}^{n-r}|\pi|\esp{\frac1{h_k^2}|\EFp{\tk{}}{\errYk{}}|^2 +|\EFp{\tk{}}{\errZk{}}|^2}
		+ 
\sum_{k=n-r+1}^{n}\esp{|\dY_{k}|^2+|\pi||\dZk{}|^2} 
\Big)
\end{align}

2.a Let us  define
\begin{align*}
\delta_i &:= \sum_{j=i}^{n-r}\esp{|\dYj{}|^2} \;,
\\
\theta_i &:=\sum_{k=i}^{n-r}|\pi|\esp{\frac1{h_k^2}|{\errYk{}}|^2 +|{\errZk{}}|^2}
		+ 
\sum_{k=n-r+1}^{n}\esp{|\dY_{k}|^2+ |\pi| |\dZk{}|^2} \;.
\end{align*}
Equation \eqref{eq stab Y step 1.c-1} reads then
\begin{align}\label{eq stab Y step 2-1}
\delta_i-\delta_{i+1} \le C|\pi|\delta_{i} + C\theta_i
\end{align}
Using a discrete version of Gronwall Lemma, we then compute
\begin{align*}
\delta_i \le C(\delta_{n-r} + \sum_{k=i}^{n-r} \theta_k e^{C(n-r-k)|\pi|})
\end{align*}
Since $\delta_{n-r} \le \eta_i$ and $\theta_k \le \theta_i$ for $k \ge i$, we compute
\begin{align*}
\delta_i \le C \theta_i \frac{1}{1-e^{C|\pi|}}
\end{align*}
This last equation combined with \eqref{eq stab Y step 2-1} leads to
\begin{align*}
\esp{|\dYi{}|^2} \le C \theta_i
\end{align*}
which concludes the proof for the $Y$-part.

2.b For the $Z$-part, the proof is concluded pluging last inequality in \eqref{eq stab interm Z}, with $i=0$ in this equation.
\eproof

\begin{Remark} It is easily checked that $\HYP{c}$ implies that
the roots of the following polynomial equations
\begin{align*}
y^{r+1}-\sum_{j=1}^r a_j y^{r-j+1} = 0\,.
\end{align*}
are in the closed unit disc and the multiple roots are in the open unit disc.
\\
It is known that in the ODEs framework this is a necessary and sufficient condition to get stability of linear multi-step schemes, see e.g \cite{but08, dem06}.
In our context, this condition is only necessary. We have to imposed $\HYP{c}$ essentially because we need to deal with the new process $Z$.
\end{Remark}

\begin{Remark} \label{re pr L2 stab}
Proposition \ref{pr L2 stab} is generic in the sense that we do not use  the particular property of the probability space nor the fact that $(\cF_{t})_{t \in [0,T]}$ is a Brownian filtration. We will use this property in the last section of this paper. 
\end{Remark}




\subsection{Study of the order}

\subsubsection{Definitions}

To study the order of the schemes, we use the following  definition of truncation errors.

The \emph{local truncation error} for the pair $(Y,Z)$ defined as 
\begin{align}
\eta_i := \eta^Y_i + \eta^Z_i,\ \ \ \ 
 (\eta^Y_i ,\eta^Z_i):=\left(\frac1{\hi^2} \esp{|Y_{\ti{}} - \check{Y}_{\ti{}}|^2},\esp{ |Z_{\ti{}} -  \check{Z}_{\ti{}}|^2}\right) ,\; i\le n-r\;, \label{eq de trunc error loc YnZ}
\end{align}
with
\begin{align}
\check{Y}_{\ti{}} &=\EFp{\ti{}}{\sum_{j=1}^r a_jY_{\ti{+j}} + h_i \sum_{j=1}^r b_{i,j} f(Y_{\ti{+j}},Z_{\ti{+j}}) + h_i b_{i,0}f(\check{Y}_{\ti{}},\check{Z}_{\ti{}}) } \label{eq def hY}
\\
\check{Z}_{\ti{}} &= \EFp{\ti{}}{ 
		\sum_{j = 1}^r \alpha_j H^\psi_{\ti{},jh} Y_{\ti{+j}}  + h_j \sum_{j=1}^r \beta_{i,j} H^\phi_{\ti{},jh} f(Y_{\ti{+j}},Z_{\ti{+j}}) \label{eq def hZ}
}
\end{align}
where $\psi$, $\phi$ belongs to $\cB^0$.

\vspace{1mm}
The \emph{global truncation error} for a given grid $\pi$ is given by
\begin{align}
 \cT(\pi) := \cT_Y(\pi) +  \mathcal{T}_Z(\pi),\ \ \ 
 (\cT_Y(\pi),  \cT_Z(\pi)) :=\left( 
\sum_{i=0}^{n-r} \hi \eta^Y_i,\ \ \  
\sum_{i=0}^{n-r} \hi \eta^Z_i \right),
 \label{eq de trunc error YnZ}
\end{align}
where $\cT_Y$ is the global truncation error for $Y$ and  $\cT_Z$ is the global truncation error for $Z$ defined as above.
%
%
%
%
%
%
%
%
%
%

\begin{Definition}\label{de order}
An approximation is said to have a \emph{global truncation error} of order $m$ if we have 
\begin{align*}
 \cT(\pi)  \le C |\pi|^{2m}
\end{align*}
for all sufficiently smooth\footnote{The required regularity assumptions will be stated in the Theorems below.} solutions to  \eqref{b_pde_t} and 
all partitions $\pi$ with sufficiently small mesh size.
\end{Definition}

\subsubsection{Expansion of the truncation error}

We study the order of the methods given in Definition \ref{de multi-step scheme} using It\^o-Taylor expansions  \cite{klopla92}. This requires the smoothness of the value function $u$ introduced in \eqref{b_pde_t}-\eqref{b_pde_T}. In order to state precisely these assumptions, we recall some notations of Chapter 5 (see  Section 5.4) in \cite{klopla92}. 

\vspace{1mm}
Let
\begin{align*}
\cM := \set{\oslash}\cup \bigcup_{m =1}^\infty\set{0,\dots,d}^m
\end{align*}
be the set of multi-indices with entries in $\set{0,\dots,d}$ endowed with the measure $\ell$ of the length of a multi-index ($\ell(\oslash)=0$ by convention).

We  introduce the concatenation operator $*$ on $\cM$ for  multi-indices with finite length: $\alpha= (\alpha_1,\dots,\alpha_p)$,  $\beta= (\beta_1,\dots,\beta_q)$ then $\alpha*\beta = (\alpha_1,\dots,\alpha_p,\beta_1,\dots,\beta_q)$. 


\vspace{1mm}
A non empty subset $\cA \subset \cM$ is called a hierarchical set if
\begin{align*}
\sup_{\alpha} \ell(\alpha) < \infty \text{ and } -\alpha \in \cA , \; \forall \alpha \in \cA \setminus \set{\oslash}
\end{align*}

For any  hierarchical $\cA$ set, we consider the remainder set $\cB(\cA)$ given by
\begin{align*}
\cB (\cA) := \set{\alpha \in \cM\setminus \cA | -\alpha \in \cA}
\end{align*}

We will use in the sequel the following sets of multi-indices, for $n\ge 0$:
\begin{align*}
\cA_{n} := \set{ \alpha \; | \; \ell(\alpha) \le n}
\end{align*}
and observe that $\cB(\cA_{n})=\cA_{n+1}\setminus \cA_{n}$.

\vspace{1mm}
For $j \in \set{1, \dots,d}$, we consider the operators:
\begin{align*}
L^{(j)} = \sum_{k=1}^d \sigma^{kj} \partial_{x_k}.
\end{align*} 

For a multi-index $\alpha=(\alpha_1,\dots,\alpha_p)$, the iteration of these operators has to be understood in the following sense
\begin{align*}
L^\alpha := L^{(\alpha_1)}\circ \dots \circ L^{(\alpha_p)}.
\end{align*}
By convention, $L^{\oslash}$ is the identity operator, recall also the definition of $L^{(0)}$ given in \eqref{operator}.
One can observe that $L^{\alpha*\beta}=L^\alpha \circ L^\beta$.

For a multi-index with finite length $\alpha$, we consider the set $\cG^\alpha$ of function $v:[0,T]\times \R^d \rightarrow \R$ for which $L^\alpha v$ is well defined and continuous. We also introduce $\cG^\alpha_b$ the subset of function $v \in \cG^\alpha$ such that the function $L^\alpha v$ is bounded.  For $v \in \cG^\alpha$, we denote  $L^\alpha u$ by $u^\alpha$.

Finally, for $n\ge 1$, we define the set $\cG^n_b$ of function $u$ such that $u^\alpha \in \cG^\alpha_b$ for all $\alpha \in \cA_n\setminus \set{\oslash}$.

\vspace{2mm}

The two following Propositions are key results to prove the high order rate of convergence of the schemes. We refer to \cite{chacri12} for proofs.

\begin{Proposition}\label{pr weak expansion Y}
Assume $d=1$. Let $m\ge 0$, then for a function $v \in \cG^{m+1}_b$,
 we have that
\begin{align*}
\EFp{t}{v(t+h,X_{t+h})} 
&=      v_t + hv_t^{(0)} + \frac{h^2}{2}v_t^{(0,0)} + \dots + \frac{h^m}{m!}v_t^{(0)_m}
                        +
        O_t(h^{m+1})
\end{align*}
\end{Proposition}



\vspace{2mm}
\begin{Proposition}\label{pr weak expansion Z} Assume $d=1$.
(i) Let $m \ge 0$, for $\psi \in \cB^m_{[0,1]}$, assuming that $v \in \cG^{m+2}_b$, 
we have 
\begin{align*}
\EFp{t}{H^{\psi}_{t, h}v(t+ h,X_{t+ h})} 
&= v^{(1)}_t + h v^{(1,0)}_t + \dots 
        + \frac{h^{m}}{m!}v^{(1)*(0)_{m}}_t  + O_t(h^{m+1})
\end{align*}
\\
(ii) For $\psi \in \cB^0_{[0,1]}$, assuming that $v \in \cG^1_b$, we have
\begin{align*}
\EFp{t}{H^{\psi}_{t, h}v(t+ h,X_{t+ h})} 
&= O_t(1)\,.
\end{align*}
\\
(iii) If $L^{(0)}\circ L^{(1)} = L^{(1)}\circ L^{(0)}$, then the expansion of (i) holds true with $\psi = 1$.
\end{Proposition}

\subsubsection{Sufficient condition for Order $m$}
\label{subsubse order m}

For the reader's convenience, we assume in this paragraph a constant time step for the grid $\pi$ i.e. $\hi = h = |\pi| := \frac{T}{n}$, for all $i$ and that the coefficients $b$, $\beta$ do not depend of $i$.

Under these conditions, the scheme given in Definition \ref{de multi-step scheme}, recalling \eqref{eq def H scheme}, rewrites, for $i\le n-r$,
\begin{align} \label{eq scheme general}
\left \{
\begin{array}{rcl}
\Yi{} &=& \EFp{\ti{}}{ \sum_{j=1}^r a_{j}\Yi{+j} + h \sum_{j=0}^r b_{j} f(\Yi{+j},\Zi{+j}) }
\vspace{2mm}
\\
\Zi{} &=&\EFp{\ti{}}{ 
		\sum_{j = 1}^r \alpha_{j} H^\psi_{\ti{},jh}  \Yi{+j} + h \sum_{j=1}^r \beta_{j} H^\phi_{\ti{},jh}  f(\Yi{+j},\Zi{+j}) 
}
\end{array}
\right.
\end{align}


%
%

\vspace{2mm}

\begin{Proposition} (Order m) \label{pr order m}
For $m\ge 2$, assume that the following  holds
\begin{align*}
(C^Y)_m&:\;\; \sum_{j=1}^r a_j j^p - p\sum_{j=0}^r b_j j^{p-1} = 0,\; 1\le p \le m
\\
 \text{ and } \quad (C^Z)_m& :\;\; \sum_{j=1}^r \alpha_j j^p - p \beta_j j^{p-1} = 0,\; 1\le p \le m-1
\end{align*} and that $u \in \cG^{m+1}_b$, then we have
\begin{align*}
\cT_Y(\pi) + \cT_Z(\pi) \le C |\pi|^{2m},
\end{align*}
provided that $\psi \in \cB^{m-1}_{[0,1]}$ and $\phi \in \cB^{m-2}_{[0,1]}$, recalling \eqref{eq scheme general}.
\end{Proposition}

\proof 

1. We first study the truncation error for the Z-part. We have that
\begin{align*}
\check{Z}_{\ti{}} &=\EFp{\ti{}}{ 
		\sum_{j = 1}^r \alpha_j H^\psi_{\ti{},jh}  Y_{\ti{+j}}  + h \sum_{j=1}^r \beta_j H^\phi_{\ti{},jh}  f(Y_{\ti{+j}},Z_{\ti{+j}})
		}
\\
&= \EFp{\ti{}}{\sum_{j=1}^r \alpha_j H^\psi_{\ti{},jh}   u(\ti{+j},X_{\ti{+j}}) - h \sum_{j=1}^r \beta_j H^\phi_{\ti{},jh}  u^{(0)}(\ti{+j},X_{\ti{+j}})} \;.
\end{align*}

Using Proposition \ref{pr weak expansion Z}, we compute
\begin{align*}
\check{Z}_{\ti{}}  = \sum_{p=0}^{m-1} \sum_{j=1}^r \alpha_j j^p \frac{h^p}{p!} u^{(1)*(0)_p}(\ti{},X_{\ti{}})
-  \sum_{p=0}^{m-2} \sum_{j=1}^r \beta_j  j^{p} \frac{h^{p+1}}{p!} u^{(1)*(0)_{p+1}}(\ti{},X_{\ti{}}) +O_{\ti{}}(|\pi|^{m})
\end{align*}
which leads to
\begin{align*}
\check{Z}_{\ti{}} & - Z_{\ti{}} = (\sum_{j=1}^r \alpha_j-1) u^{(1)}(\ti{},X_{\ti{}})
+ \sum_{p=1}^{m-1} \frac{h^p}{p!}u^{(1)*(0)_p}(\ti{},X_{\ti{}})(\sum_{j=1}^r \alpha_j j^p - p\sum_{j=1}^r j^{p-1}\beta_j )
+O_{\ti{}}(|\pi|^{m})
\end{align*}
Under $(C^Z)_m$, we obtain
\begin{align*}
\check{Z}_{\ti{}} & - Z_{\ti{}} = O_{\ti{}}(|\pi|^{m}) 
\end{align*}
which leads directly to
\begin{align}
\eta^Z_i = O(|\pi|^{2m}),\quad i\le n-r. \label{eq trunc error gene Z}
\end{align}

2.a  We now study the truncation error for the Y-part. Let us introduce 
\begin{align*}
\bar{Y}_{\ti{}} &=\EFp{\ti{}}{\sum_{j=1}^r a_jY_{\ti{+j}} + h \sum_{j=0}^r b_{j} f(Y_{\ti{+j}},Z_{\ti{+j}}) }
\end{align*}
We have that
\begin{align*}
\check{Y}_{\ti{}} &= \bar{Y}_{\ti{}} +  h b_{0} \Big(f(\check{Y}_{\ti{}},\check{Z}_{\ti{}}) - f(Y_{\ti{}},Z_{\ti{}})\Big)
\end{align*}
Since $f$ is Lipschitz-continuous, we get that for $|\pi|$ small enough,
\begin{align} \label{eq majo Y 1}
\check{Y}_{\ti{}} - Y_{\ti{}} = O_{\ti{}}(\bar{Y}_{\ti{}} - Y_{\ti{}}) + |\pi| O_{\ti{}}(\check{Z}_{\ti{}} - Z_{\ti{}}) \;.
\end{align}

2.b Now observe that
\begin{align*}
\bar{Y}_{\ti{}} &= \EFp{\ti{}}{\sum_{j=1}^r a_j u(\ti{+j},X_{\ti{+j}}) - h \sum_{j=0}^r b_j u^{(0)}(\ti{+j},X_{\ti{+j}}) } \;.
\end{align*}

Using Proposition \ref{pr weak expansion Y}, we compute
\begin{align*}
\bar{Y}_{\ti{}}  = \sum_{p=0}^m \sum_{j=1}^r a_j j^p \frac{h^p}{p!} u^{(0)_p}(\ti{},X_{\ti{}})
-  \sum_{p=0}^{m-1} \sum_{j=0}^r b_j  j^{p} \frac{h^{p+1}}{p!} u^{(0)_{p+1}}(\ti{},X_{\ti{}}) +O_{\ti{}}(|\pi|^{m+1})
\end{align*}
which leads to
\begin{align*}
\bar{Y}_{\ti{}} & - Y_{\ti{}} = (\sum_{j=1}^r a_j-1) u(\ti{},X_{\ti{}})
+ \sum_{p=1}^m \frac{h^p}{p!}u^{(0)_p}(\ti{},X_{\ti{}})(\sum_{j=1}^r a_j j^p - p\sum_{j=0}^r j^{p-1}b_j )
+O_{\ti{}}(|\pi|^{m+1}).
\end{align*}
Under $(C^Y)_m$, we thus get 
\begin{align*}
\bar{Y}_{\ti{}} & - Y_{\ti{}} = O_{\ti{}}(|\pi|^{m+1}) 
\end{align*}

2.c Combining the last inequality with \eqref{eq majo Y 1} and \eqref{eq trunc error gene Z}, we then obtain
\begin{align*}
\eta^Y_i = O(|\pi|^{2m+2})\;,\; i \le n-r\;.
\end{align*}

\vspace{2mm}
3. Combining the last equation with \eqref{eq trunc error gene Z}, we conclude that 
\begin{align*}
\cT(\pi)= O(|\pi|^{2m}) ,
\end{align*}
and so the scheme is of order $m$.

\eproof

\subsection{Convergence results and examples of high order methods}
\begin{Theorem} \label{th main conv res} Under \HYP{c}, assuming that the scheme is of order $m$ according to Definition \ref{de order} and that
\begin{align}
\max_{0\le j \le r-1} \esp{|Y_{t_{n-j}} - Y_{n-j}|^2 + h|Z_{t_{n-j}}- Z_{n-j}|^2} \le C|\pi|^{2m} \label{eq order cond term}
\end{align}
we have
\begin{align*}
\cE_Y(\pi) + \cE_Z(\pi) \le C |\pi|^{2m}\;. 
\end{align*}
\end{Theorem}

\vspace{2mm}
\proof
We simply observe that the solution $(Y,Z)$ of the BSDE is also the solution of a perturbed scheme with $\errYi{}:= \check{Y}_{\ti{}}-Y_{\ti{}}$ and $\errZi{}:= \check{Z}_{\ti{}}-Z_{\ti{}}$. The proof then follows directly from Proposition \ref{pr L2 stab}.
\eproof

\vspace{6mm}
In particular, in the special setting of paragraph \ref{subsubse order m}, a straightforward application of Theorem \ref{th main conv res} and Proposition \ref{pr order m} leads to
\begin{Corollary}\label{co gen res}
Under \HYP{c} and $(C^Y)_m$- $(C^Z)_m$,   assuming that  \eqref{eq order cond term} holds, we have
\begin{align*}
\cE_Y(\pi) + \cE_Z(\pi) \le C |\pi|^{2m}\;, 
\end{align*}
provided $u \in \cG^{m+1}_b$ and $\psi \in \cB^{m-1}_{[0,1]}$, $\phi \in  \cB^{m-2}_{[0,1]}$.
\end{Corollary}

\vspace{6mm}
To illustrate the previous results, we conclude this section by giving two examples of high order method which can be designed using Corollary \ref{co gen res}. 
\begin{Example} (Nystrom's method) \label{ex leapfrog}
The following scheme is --for the Y-part-- inspired by the 
\emph{Leap-frog (or Nystrom's) method} for ODE, namely
\begin{align*} 
\left \{
\begin{array}{rcl}
\Yi{} &=& \EFp{\ti{}}{ \Yi{+2} + 2h f(\Yi{+1},\Zi{+1}) }
\vspace{2mm}
\\
\Zi{} &=&\EFp{\ti{}}{ 
		H^\psi_{\ti{},2h}  \Yi{+2} + 2 h H^\phi_{\ti{},2h}  f(\Yi{+2},\Zi{+2})}.
\end{array}
\right.
\end{align*}
This 2-step method  is convergent and the rate of convergence  is at least of order 2, assuming that
$u \in \cG^3_b$ and $\psi \in \cB^1_{[0,1]}$, $\phi \in \cB^0_{[0,1]}$.
\end{Example}

\vspace{4mm}
\begin{Example} (Milne's method) \label{ex milne's}
The second scheme we propose is inspired --for the Y-part-- by the \emph{Milne's method} for ODE, namely
\begin{align*} 
\left \{
\begin{array}{rcl}
\Yi{} &=& \EFp{\ti{}}{ \Yi{+4} + h\Big(\frac83 f(\Yi{+1},\Zi{+1}) - \frac43f(\Yi{+2},\Zi{+2}) + \frac83 f(\Yi{+3},\Zi{+3}) \Big) }
\vspace{2mm}
\\
\Zi{} &=&
\EFp{\ti{}}{ 
		H^\psi_{\ti{},4h} \Yi{+4}  +  h \Big( \frac83H^\phi_{\ti{},h} f(\Yi{+1},\Zi{+1})  -\frac43H^\phi_{\ti{},2h} f(\Yi{+2},\Zi{+2}) + \frac83H^\phi_{\ti{},3h} f(\Yi{+3},\Zi{+3}) \Big) }
. 
\end{array}
\right.
\end{align*}
This 4-step method is convergent and  the rate of convergence is at least of order 4,
assuming that
$u \in \cG^5_b$ and $\psi \in \cB^3_{[0,1]}$, $\phi \in \cB^2_{[0,1]}$.
\end{Example}

%
%


\section{Adams Methods}

In this section, we introduce methods for BSDEs
inspired by Adams methods from the ODE framework. These methods are of two kinds: explicit methods , also called Adams-Bashforth, or implicit methods, also called, Adams-Moulton.

The schemes introduced in Definition \ref{de multi-step scheme} are always explicit for the $Z$-part but may be implicit for the $Y$-part. So, for the $Z$-part, we  use Adams-Bashforth approximation which may then be combined with explicit or implicit approximation for the $Y$-part.

\vspace{2mm}
We first study methods combining Adams-Moulton type approximation for the $Y$-part and  Adams-Bashforth type approximation for the $Z$-part. We  show that these methods are really efficient because high order rate of convergence can be achieved, assuming smoothness of the value function.
We then quickly discuss the case of explicit methods, i.e. Adams-Bashforth type approximation both for the $Y$-part and $Z$-part. 

At the end of this section, we use these Adams type approximation to design Predictor-Corrector methods for BSDEs.

\subsection{Implicit methods}

These methods are inspired by Adams-Moulton method  for the $Y$-part and Adams-Bashforth for the $Z$-part. 

They have the following form, for $i\le n-r$,
\begin{align*} 
(AMB)_r\,: \;\left \{
 \begin{array}{rcl}
\Yi{} &=& \EFp{\ti{}}{ \Yi{+1} + h_i \sum_{j=0}^r b_{i,j,r} f(\Yi{+j},\Zi{+j}) }
\vspace{2mm}
\\
\Zi{} &=&\EFp{\ti{}}{ 
		H^\psi_{\ti{},h}  \Yi{+1 } + h_i \sum_{j=1}^r \beta_{i,j,r} H^\phi_{\ti{},jh} f(\Yi{+j},\Zi{+j}) 
}
\end{array}
\right.
\end{align*}
where $\psi,\phi \in \cB^0_{[0,1]}$.

\vspace{2mm}
The coefficients  for the $Y$-part are given by
\begin{align}
b_{i,j,r} = \frac1{\hi}\int_{\ti{}}^{\ti{+1}} L_{i,j,r}(s) \ud s\, ,\; \text{ with } L_{i,j,r}(t) = \Prod_{k=0,k \neq j}^{r} \frac{t-\ti{+k} }{\ti{+j}-\ti{+k} } \;, \; 0 \le j \le r. \label{eq de b AM}
\end{align}
The Lagrange polynomials $L_{i,j,r}$ are  of degree $r$ and $L_{i,j,r}(\ti{+j})=1$, which implies
\begin{align}
\sum_{j=0}^r (\ti{+j}-\ti{})^k L_{i,j,r}(t) = (t-\ti{})^k\,,\quad 0\le k \le r\,. \label{eq useful Lijr}
\end{align}
The definition of the $b$-coefficients means that
\begin{align*}
\Yi{} = \EFp{\ti{}}{ \Yi{+1} + \int_{\ti{}}^{\ti{+1}} Q^Y_{i,r}(t) \ud t }
\end{align*}
where $Q^Y_{i,r}$ is a polynomial of degree less than $r$ satisfying
\begin{align*}
Q^Y_{i,r}(\ti{+j}) = f(\Yi{+j},\Zi{+j})\;,\; 0 \le j \le r\;.
\end{align*}

In the case where the time step is constant, the coefficient does not depends on $i$ and are given by
\begin{align*}
b_{j,r} = \int_0^1 \ell_{j,r}(s) \ud s\, , \; \text{ with } \ell_{j,r}(s) = \Prod_{k=0,k \neq j}^{r+1} \frac{s-k}{j+1-k} \;,\; 0 \le j \le r\;.
\end{align*}

The coefficients  for the $Z$-part are given by
\begin{align}
\beta_{i,j,r} = \frac1{\hi}\int_{\ti{}}^{\ti{+1}} \tilde{L}_{i,j,r}(s) \ud s\, , \; 1 \le j \le r, \text{ with } \tilde{L}_{i,j,r}(s) = \Prod_{k=1,k \neq j}^r \frac{t-\ti{+k} }{\ti{+j}-\ti{+k} } \;. \label{eq de tilde L}
\end{align}
The Lagrange polynomials  $\tilde{L}_{i,j,r}$ are  of degree $r-1$ and $\tilde{L}_{i,j,r}(\ti{+j})=1$, which implies
\begin{align}
\sum_{j=1}^r (\ti{+j}-\ti{})^k \tilde{L}_{i,j,r}(t) = (t-\ti{})^k\,,\quad 0\le k \le r-1\,. \label{eq useful tilde Lijr}
\end{align}

The definition of the $\beta$-coefficients means that
\begin{align*}
\Zi{} = \EFp{\ti{}}{ H^\psi_{\ti{},h} \Yi{+1} + \int_{\ti{}}^{\ti{+1}} Q^Z_{i,r}(t) \ud t }
\end{align*}
where $Q^Z_{i,r}$ is a polynomial of degree less than $r-1$ satisfying
\begin{align*}
Q^Z_{i,r}(\ti{+j}) = H^\phi_{\ti{},jh} f(\Yi{+j},\Zi{+j})\;,\; 1 \le j \le r\;.
\end{align*}

In the case where the time step is constant, the coefficient does not depends on $i$ and are given by
\begin{align*}
\beta_{j,r} = \int_0^1 \ell_{j,r}(s) \ud s\, , \; \text{ with } \ell_{j,r}(s) = \Prod_{k=1,k \neq j}^r \frac{s-k}{j-k} \;.
\end{align*}

When the time step is constant, the table below gives the $b$-coefficients and $\beta$-coefficients for $r \le 4$:
\begin{align*}
\begin{array}{c|ccccc|cccc}
r & b_0 & b_1 & b_2 & b_3 & b_4& \beta_1 & \beta_2 & \beta_3 & \beta_4 
\\ \hline
1 & \frac12 & \frac12 & &  & & 1 & & &
\\
2 & \frac52 & \frac8{12} & - \frac1{12}& & & \frac32 & - \frac12 & &
\\
3 &  \frac9{24}& \frac{19}{24}  & -\frac5{24} & \frac1{24} & & \frac{23}{12} & - \frac{16}{12} & \frac{5}{12} &
\\
4 & \frac{251}{720} & \frac{646}{720}  &-\frac{264}{720}  & \frac{106}{720}  & -\frac{19}{720} &  \frac{55}{24} & - \frac{59}{24} & \frac{37}{24} & - \frac9{24}
\end{array}
\end{align*}

\vspace{2mm}
\begin{Proposition}\label{pr AM r}
The $(AMB)_r$ method is convergent and at least of order $r+1$, provided that $\psi \in \cB^{r}$, $\phi \in \cB^{r-1}$ and $u \in \cG^{r+2}_b$.
\end{Proposition}

\proof 1. The stability of the schemes comes from a direct application of Proposition \ref{pr L2 stab}, since obviously $\HYP{c}$ holds for $(AMB)_r$. Following Theorem \ref{th main conv res}, we only have to study the order of the method. 

2.a We first study the error for the $Z$ part. Observe that, recalling \eqref{eq def hZ},
\begin{align*}
\check{Z}_{\ti{}} &:= \EFp{\ti{}}{H^\psi_{\ti{},h} Y_{\ti{+1}} + \sum_{j=1}^r H^\phi_{\ti{},jh}  f(Y_{\ti{+j}},Z_{\ti{+j}})\int_{\ti{}}^{\ti{+1}} L_{i,j,r}(t) \ud t }
\\
&= \EFp{\ti{}}{H^\psi_{\ti{},h}  u(\ti{+1},X_{\ti{+1}}) - \sum_{j=1}^r H^\phi_{\ti{},jh} u^{(0)}(\ti{+j},X_{\ti{+j}})\int_{\ti{}}^{\ti{+1}} L_{i,j,r}(t) \ud t }
\end{align*}
Using Proposition \ref{pr weak expansion Z}, we get
\begin{align*}
\check{Z}_{\ti{}} - Z_{\ti{}} = \sum_{k=1}^{r}\frac{\hi^k}{k!}u^{(1)*(0)_k}(\ti{},X_{\ti{}}) - \sum_{j=1}^r  \int_{\ti{}}^{\ti{+1}} \tilde{L}_{i,j,r}(t) \ud t \sum_{k=0}^{r-1} \frac{u^{(1)*(0)_{k+1}}(\ti{},X_{\ti{}})}{k!}(\ti{+j}-\ti{})^{k}  + O_{\ti{}}(|\pi|^{r+1})
\end{align*}
which reads also
\begin{align*}
\check{Z}_{\ti{}} - Z_{\ti{}} = \sum_{k=1}^{r}\frac{\hi^k}{k!}u^{(1)*(0)_k}(\ti{},X_{\ti{}}) - \sum_{k=0}^{r-1} \frac{u^{(1)*(0)_{k+1}}(\ti{},X_{\ti{}})}{k!}  \int_{\ti{}}^{\ti{+1}} \sum_{j=1}^r (\ti{+j}-\ti{})^{k}  \tilde{L}_{i,j,r}(t) \ud t   + O_{\ti{}}(|\pi|^{r+1})\;.
\end{align*}
Using \eqref{eq useful tilde Lijr}, 
we obtain
\begin{align*}
\check{Z}_{\ti{}} - Z_{\ti{}} & = \sum_{k=1}^r\Big(\frac{\hi^k}{k!}- \frac1{(k-1)!}\int_{\ti{}}^{\ti{+1}} (t-\ti{})^{k-1}\ud t \Big)u^{(1)*(0)_k}(\ti{},X_{\ti{}})  + O_{\ti{}}(|\pi|^{r+1})
\\
& = O_{\ti{}}(|\pi|^{r+1}).
\end{align*}

2.b We now study the truncation error for the $Y$ part. 
Let us define,
\begin{align*}
\bar{Y}_{\ti{}} &:= \EFp{\ti{}}{Y_{\ti{+1}} + \sum_{j=0}^r f(Y_{\ti{+j}},Z_{\ti{+j}})\int_{\ti{}}^{\ti{+1}} L_{i,j,r}(t) \ud t }
\end{align*}

Observe that 
\begin{align*}
\check{Y}_{\ti{}} &:= \bar{Y}_{\ti{}}   + \Big( f(\check{Y}_{\ti{}},\check{Z}_{\ti{}})-f(Y_{\ti{}},Z_{\ti{}}) \Big) \int_{\ti{}}^{\ti{+1}} L_{i,0,r}(t) \ud t
\end{align*}
which leads since $f$ is Lipschitz continuous, for $|\pi|$ small enough, to
\begin{align}\label{eq adams interm 1}
 \check{Y}_{\ti{}}  - Y_{\ti{}}= O_{\ti{}}(|\bar{Y}_{\ti{}}- Y_{\ti{}}|) + |\pi|O_{\ti{}}(|\check{Z}_{\ti{}}- Z_{\ti{}}|)\;.
\end{align}

Now,
\begin{align*}
\bar{Y}_{\ti{}}  &= \EFp{\ti{}}{u(\ti{+1},X_{\ti{+1}}) - \sum_{j=0}^r u^{(0)}(\ti{+1},X_{\ti{+1}})\int_{\ti{}}^{\ti{+1}} L_{i,j,r}(t) \ud t }
\end{align*}
Using Proposition \ref{pr weak expansion Y}, we get
\begin{align*}
\bar{Y}_{\ti{}} - Y_{\ti{}} = \sum_{k=1}^{r+1}\frac{\hi^k}{k!}u^{(0)_k}(\ti{},X_{\ti{}}) - \sum_{j=0}^r  \int_{\ti{}}^{\ti{+1}} L_{i,j,r}(t) \ud t \sum_{k=0}^{r} \frac{u^{(0)_{k+1}}(\ti{},X_{\ti{}})}{k!}(\ti{+j}-\ti{})^{k}  + O_{\ti{}}(|\pi|^{r+2})
\end{align*}
which reads also
\begin{align*}
\bar{Y}_{\ti{}} - Y_{\ti{}} = \sum_{k=1}^{r+1}\frac{\hi^k}{k!}u^{(0)_k}(\ti{},X_{\ti{}}) - \sum_{k=0}^{r} \frac{u^{(0)_{k+1}}(\ti{},X_{\ti{}})}{k!}  \int_{\ti{}}^{\ti{+1}} \sum_{j=0}^r (\ti{+j}-\ti{})^{k}  L_{i,j,r}(t) \ud t   + O_{\ti{}}(|\pi|^{r+2})\;.
\end{align*}
Using \eqref{eq useful  Lijr}, 
we obtain
\begin{align*}
\bar{Y}_{\ti{}} - Y_{\ti{}} & = \sum_{k=1}^{r+1}\Big(\frac{\hi^k}{k!}- \frac1{(k-1)!}\int_{\ti{}}^{\ti{+1}} (t-\ti{})^{k-1}\ud t \Big)u^{(0)_k}(\ti{},X_{\ti{}})  + O_{\ti{}}(|\pi|^{r+2})
\\
& = O_{\ti{}}(|\pi|^{r+2}).
\end{align*}
Thus, using \eqref{eq adams interm 1},
\begin{align*}
\check{Y}_{\ti{}} - Y_{\ti{}}  = O_{\ti{}}(|\pi|^{r+2}).
\end{align*}

2.c Combining the results of steps 1.a and 1.b, we obtain
\begin{align*}
\cT(\pi)= O_{\ti{}}(|\pi|^{r+1})\;.
\end{align*}
which concludes the proof.
\eproof

\subsection{Explicit methods}
These methods are inspired by Adams-Bashforth method both  for the $Y$-part and $Z$-part.  
\begin{align*} 
(ABB)_r \left \{
\begin{array}{rcl}
\Yi{} &=& \EFp{\ti{}}{ \Yi{+1} + \hi \sum_{j=1}^r b_{i,j,r} f(\Yi{+j},\Zi{+j}) }
\vspace{2mm}
\\
\Zi{} &=&\EFp{\ti{}}{ 
		H^\psi_{\ti{},h}  \Yi{+1 } + h_i \sum_{j=1}^r \beta_{i,j,r} H^\phi_{\ti{},jh} f(\Yi{+j},\Zi{+j}) 
}
\end{array}
\right.
\end{align*}
where $\psi, \phi \in \cB^0_{[0,1]}$.

Now, the coefficients  for the $Y$-part are  given by
\begin{align*}
b_{i,j,r} = \frac1{\hi}\int_{\ti{}}^{\ti{+1}} \tilde{L}_{i,j,r}(s) \ud s\, ,
\end{align*}
recalling \eqref{eq de tilde L}.




\vspace{2mm}
From the proof of Proposition \ref{pr AM r}, step 1.b. we know that we can obtain a truncation error for the  $Z$-part s.t. $\hat{Z_{\ti{}}} - Z_{\ti{}} =O_{\ti{}}(|\pi|^{r+1})$. But here, due to the explicit feature of the $Y$ part and thus an order $r$ global error only, we only need to retrieve an error for the $Z$-part of order $r$ as well. This simply means that the scheme, for the $Z$-part, has one more coefficient than needed. So one can set 
\begin{align*}
\beta_{i,j,r} = b_{i,j,r}
\end{align*}
or 
\begin{align*}
\beta_{i,r,r}= 0 \text{ and } \beta_{i,j,r} = b_{i,j,r-1}, \; 1\le j \le r-1\;.
\end{align*}

Following the arguments of the proof of Proposition \ref{pr AM r}, one obtains
\begin{Proposition}\label{pr AB r} The $(ABB)_r$ method is convergent and at least of order $r$, provided that $\psi \in \cB^{r-1}$, $\phi \in \cB^{r-2}$ and $u \in \cG^{r+1}_b$.
\end{Proposition}



\subsection{Predictor-Corrector methods}

These methods are fully explicit method but have a better rate of convergence than the $(ABB)_r$ methods presented above. Nevertheless, they require the computation of one more conditional expectation by step. This has to be compared  in practice to the Picard Iteration required by $(AMB)_r$ approximation.
\begin{Definition} \label{de pece}
\begin{align} 
(PC)_r \left \{
\begin{array}{rcl}
\Zi{} &=&\EFp{\ti{}}{ 
		H^\psi_{\ti{},h}  \Yi{+1 } + \hi \sum_{j=1}^r \beta_{i,j,r} H^\phi_{\ti{},jh}  f(\Yi{+j},\Zi{+j}) }
\vspace{2mm}
\\
{}^p\Yi{} &=& \EFp{\ti{}}{ \Yi{+1} + \hi \sum_{j=1}^r \beta_{i,j,r} f(\Yi{+j},\Zi{+j}) }
\vspace{2mm}
\\
\Yi{} &=& \EFp{\ti{}}{ \Yi{+1} + \hi \sum_{j=1}^r b_{i,j,r} f(\Yi{+j},\Zi{+j}) } + \hi b_{i,0,r} f({}^p\Yi{},\Zi{}) 
\end{array}
\right. \label{eq de pece}
\end{align}
where $\psi, \phi \in \cB^0_{[0,1]}$. 
The $b$-coefficients are given by \eqref{eq de b AM} and the $\beta$-coefficients by \eqref{eq de tilde L}.
\end{Definition}

\vspace{2mm}
\begin{Theorem} \label{th conv pece}
The $(PC)_r$ method is convergent and at least of order $r+1$, provided that $\psi \in \cB^{r}$, $\phi \in \cB^{r-1}$ and $u \in \cG^{r+2}_b$.
\end{Theorem}

\vspace{2mm} As usual, the proof of this Theorem is splitted in two steps below. We first study the stability of the above schemes and then their truncation errors.

\subsubsection{Stability}

To study the \emph{stability} of the methods \eqref{eq de pece}, we introduce first a pertubed version of the scheme

\begin{align*} 
\left \{
\begin{array}{rcl}
\tZi{} &=&\EFp{\ti{}}{ 
		H^\psi_{\ti{},h}  \tYi{+1 } + \hi \sum_{j=1}^r \beta_{i,j,r} H^\phi_{\ti{},jh}  f(\tYi{+j},\tZi{+j}) + \errZi{}}
\vspace{2mm}
\\
{}^p\tYi{} &=& \EFp{\ti{}}{ \tYi{+1} + \hi \sum_{j=1}^r b_{i,j,r} f(\tYi{+j},\tZi{+j}) }
\vspace{2mm}
\\
\tYi{} &=& \EFp{\ti{}}{ \tYi{+1} + \hi \sum_{j=1}^r b^*_{i,j,r} f(\tYi{+j},\tZi{+j}) + \hi b^*_{i,0,r} f({}^p\tYi{},\tZi{}) + \errYi{}}  
\end{array}
\right.
\end{align*}

where $ \errYi{}$, $\errZi{}$ are random variables belonging to $L^2(\cF_{\ti{}})$, for $i \le n-r$.


\begin{Proposition} \label{pr pc stab}
The scheme given in \eqref{eq de pece}  is $L^2$-stable, recalling Definition \ref{de stab}.
\end{Proposition}

\proof 
For $|\pi|$ small enough, we compute, denoting ${}^p \dY = {}^pY - {}^p\tilde{Y} $,
\begin{align}
\esp{|{}^p\dYi{}|^2} & \le (1+ |\pi|)\esp{|\dYi{+1}|^2} + C\sum_{j=1}^rh_i \esp{|\dYi{+j}|^2 +|\dZi{+j}|^2 } \label{eq stab pece 1}
\\
\esp{|\dYi{}|^2} & \le (1+ |\pi|)\esp{|\dYi{+1}|^2} + C\Big(h_i \esp{|{}^p\dYi{}|^2}  +  \sum_{j=1}^rh_i \esp{|\dYi{+j}|^2 +|\dZi{+j}|^2 } \Big) + \frac{C}{h_i}|\errYi{}|^2 \label{eq stab pece 2}
\end{align}
Plugging \eqref{eq stab pece 1} into  \eqref{eq stab pece 2} and using the discrete version of Gronwall's Lemma,  we obtain 
\begin{align}\label{eq pdY interm 1}
\esp{|\dYi{}|^2} \le C\Big( |\pi|\sum_{k=i}^{n-r}\esp{|\dYk{}|^2} + \sum_{k=i}^{n-r}h_k\esp{|\dZk{}|^2} &+ 
\sum_{k=i}^{n-r}\frac1{h_k}\esp{|\errYk{}|^2} 
		+ \sum_{k=n-r+1}^{n}\esp{|\dY_{k}|^2} \Big) \;.
\end{align}
Using the same arguments as in step 1.b of the proof of Proposition \ref{pr L2 stab}, we retrieve that
\begin{align}
\sum_{k=i}^{n-r} h_k \esp{|\dZk{}|^2} &\le C\Big( \sum_{k=n-r+1}^{n}\esp{|\dYk{}|^2+|\pi||\dZk{}|^2} 
+  |\pi| \sum_{k=i}^{n-r}\esp{|\dYk{}|^2} + \sum_{k=i}^{n-r}\esp{\frac1{h_k}|{\errYk{}}|^2 +h_k|{\errZk{}}|^2}
\Big) \label{eq stab interm Z pece}
\end{align}
This leads, using \eqref{eq pdY interm 1},
\begin{align}
\esp{|\dYi{}|^2} \le C\Big( |\pi|\sum_{j=i}^{n-r}\esp{|\dYj{}|^2} 
+ 
\sum_{k=i}^{n-r}|\pi|\esp{\frac1{h_k^2}|\EFp{\tk{}}{\errYk{}}|^2 +|\EFp{\tk{}}{\errZk{}}|^2}
		+ 
\sum_{k=n-r+1}^{n}\esp{|\dY_{k}|^2+|\pi||\dZk{}|^2} 
\Big)
\end{align}
which corresponds to \eqref{eq stab Y step 1.c-1}.

The proof is then concluded using the same arguments as in step $2$ of the proof of Proposition \ref{pr L2 stab}.
\eproof

\subsubsection{Truncation error}

\begin{Proposition}\label{pr pc truncation error}
The scheme given in Definition \ref{de pece} is at least of order $r+1$ provided that $\psi \in \cB^{r}$, $\phi \in \cB^{r-1}$ and $u \in \cG^{r+2}_b$.
\end{Proposition}

\proof 1. The truncation error for the $Z$-part is the same that the one of the $(AMM)_r$ method. From the proof of Proposition \ref{pr AM r} step 1, we get
\begin{align*}
\check{Z}_{\ti{}} - Z_{\ti{}} = O_{\ti{}}(|\pi|^{r+1})\;.
\end{align*}

2. The study of the truncation error for the Y-part is slightly more involved. Let us define
\begin{align*}
Y^*_{\ti{}} :=  \EFp{\ti{}}{ Y_{\ti{+1}} + \hi \sum_{j=1}^r b_{i,j,r} f(Y_{\ti{+j}},Z_{\ti{+j}}) } + \hi b_{i,0,r} f(Y^*_{\ti{}},\hat{Z}_{\ti{}})
\end{align*}
Using the proof of Proposition \ref{pr AM r} step 2, we know that
\begin{align} \label{eq trunc error corr}
Y^*_{\ti{}} - Y_{\ti{}} = O_{\ti{}}(|\pi|^{r+2})
\end{align}
this quantity represents the truncation error for the Y-part of the Adams-Moulton method.

We also define
\begin{align*}
{}^p\check{Y}_{\ti{}}&= \EFp{\ti{}}{ Y_{\ti{+1}} + \hi \sum_{j=1}^r \beta_{i,j,r} f(Y_{\ti{+j}},Z_{\ti{+j}}) }
\\
\check{Y}_{\ti{}}&= \EFp{\ti{}}{ Y_{\ti{+1}} + \hi \sum_{j=1}^r b_{i,j,r} f(Y_{\ti{+j}},Z_{\ti{+j}}) } + \hi b_{i,0,r} f({}^p\check{Y}_{\ti{}},\check{Z}_{\ti{}}) 
\end{align*}
The term ${}^p\check{Y}_{\ti{}}- Y_{\ti{}}$ represents then the truncation error for the Predictor part, adapting the arguments of Proposition \ref{pr AM r} step 1, we have
\begin{align} \label{eq trunc error pred}
{}^p\check{Y}_{\ti{}}- Y_{\ti{}} = O_{\ti{}}(|\pi|^{r+1})
\end{align}  

The term $\check{Y}_{\ti{}}- Y_{\ti{}}$ is the truncation error we are interested in.

We then observe that
\begin{align*}
\check{Y}_{\ti{}} &= Y^*_{\ti{}} + \hi b_{i,0,r} \Big(  f({}^p\check{Y}_{\ti{}},\check{Z}_{\ti{}})-  f(Y^*_{\ti{}},\check{Z}_{\ti{}}) \Big)
\\
&= Y^*_{\ti{}} + \hi b_{i,0,r} \Big(  f({}^p\check{Y}_{\ti{}},\check{Z}_{\ti{}})-f(Y_{\ti{}},\check{Z}_{\ti{}}) + f(Y_{\ti{}},\check{Z}_{\ti{}})-  f(Y^*_{\ti{}},\check{Z}_{\ti{}}) \Big)
\end{align*}
Since $f$ is Lipschitz continuous, we obtain
\begin{align*}
\check{Y}_{\ti{}} - Y_{\ti{}} = O_{\ti{}}(| Y^*_{\ti{}} - Y_{\ti{}} |) + |\pi|O_{\ti{}}(|{}^p\check{Y}_{\ti{}}- Y_{\ti{}}| )
\end{align*}
Combining \eqref{eq trunc error corr}and \eqref{eq trunc error pred}, we then obtain
\begin{align*}
\check{Y}_{\ti{}} - Y_{\ti{}} = O_{\ti{}}(|\pi|^{r+2})\;.
\end{align*}
3. From step 1. and step 2. above, we obtain that
\begin{align*}
\cT(\pi) = O(|\pi|^{r+1}),
\end{align*}
which concludes the proof.
\eproof



\section{Numerical illustration}


In this part, we provide a numerical illustration for the  results presented above.
The scheme given in Definition \ref{de multi-step scheme} is still a theoretical one because in practice one has to compute the conditional expectation involved.
Many methods have been studied already in the context of BSDEs: regression methods \cite{goblem06}, quantization methods \cite{balpag03,delmen06}, Malliavin calculus methods \cite{boutou04, bouwar11} and tree based methods e.g. Cubature methods \cite{criman10}. 

To illustrate our previous results, i.e. the order of the time discretization error, we will focus on the simple case where $d=1$ and $X=W$, in the spirit of \cite{zhazha10}. Obviously, further numerical experiments are needed, specially in high dimension. Because we are looking towards high order approximation, it seems reasonable to combine the present multi-step schemes with Cubature methods \cite{criman10}. This is left for further research.

In the sequel, we will also assume that the $r$ terminal conditions are perfectly known. Generally, this won't be the case but it is not really a problem, see Remark \ref{re intro} (i). We explain below how the Brownian motion is approximated and give the expression of the numerical scheme which is implemented in practice. We show that this scheme is convergent and characterise its convergence order. The error we are dealing with is now composed of the discrete time error and the space discretization error. Finally, we provide some numerical results, where we compute the empirical convergence rate.


\subsection{Empirical schemes}
In order to define the scheme implemented in practice, we use a multinomial approximation of the Brownian motion.
Let us consider a discrete random variable $\xi$ matching the moments of a gaussian variable $G$ up to order $K$, i.e.
\begin{align*}
\esp{\xi^k} = \esp{G^k}, \quad 0 \le k \le K.
\end{align*}
In dimension 1, an efficient way to construct $\xi$ is to use quadrature formula.

On a (discrete, but big enough) probability space $(\widehat{\Omega}, \widehat{\P})$, we are then given $(\xi_i)_{1\le i \le n}$, i.i.d random variables with the same law as $\xi$ and define
\begin{align}
\Wpi{} = x + \sum_{j=1}^i \sqrt{t_{j}-t_{j-1}}\xi_j\;,\quad \forall \, \ti{} \in \pi\,.
\label{de approx W}
\end{align}

For later use, we say that $\Wpi{}$ is an \emph{order $K$ approximation} of the Brownian motion.

\vspace{2mm}
We also denote by $(\widehat{\cF})_{t \in \pi}$  the filtration generated by $(\widehat{W}^{0,x}_{t})_{t\in \pi}$ and $\whEFp{t}{\cdot}$ the related conditional expectation.

We can now define the numerical scheme which is used in practice.
\begin{Definition} (Linear multi-step) \label{de lms numerical scheme}

(i) To initialise the scheme with $r$ steps, $r\ge 1$, we set, for $0 \le j \le r-1$,  
$$(Y_{n-j},Z_{n-j})=(u(t_{n-j},\widehat{W}^{0,x}_{t_{n-j}}),\partial_x u (t_{n-j},\widehat{W}^{0,x}_{t_{n-j}})).$$  

(ii) For $i \le n-r$, the  computation of $(Y_i,Z_i)$ involves $r$ steps  and is given by
\begin{align} \label{de numerics scheme}
\left \{
\begin{array}{rcl}
\whYi{} &=& \whEFp{\ti{}}{ \sum_{j=1}^r a_{j}\whYi{+j}+ h_i \sum_{j=0}^r b_{j} f(\whYi{+j},\whZi{+j}) }
\vspace{2mm}
\\
\whZi{}  &=&\whEFp{\ti{}}{ 
		\sum_{j = 1}^r \whHi{,j} \Big( \alpha_{j} \whYi{+j} + h_i \sum_{j=1}^r \beta_{j}  f(\whYi{+j},\whZi{+j})  \Big)
}
\end{array}
\right.
\end{align}

The coefficient $H_{i,j}$ are the discrete version of coefficients given in \eqref{eq def H scheme}. From Proposition \ref{pr weak expansion Z} (iii), we observe that, in the case $X=W$, they can simply be defined as approximation of the Brownian increment, i.e.
\begin{align} \label{de whH}
 \whHi{,j} = \frac{\Wpi{+j}-\Wpi{}}{\ti{+j}-\ti{}}\;.
\end{align}
\end{Definition}

\vspace{2mm}
When implementing Predictor-Corrector methods, we use
\begin{Definition} (Predictor-Corrector) \label{de PC numerical scheme}
\\
(i) To initialise the scheme with $r$ steps, $r\ge 1$, we set, for $0 \le j \le r-1$,  
$$(Y_{n-j},Z_{n-j})=(u(t_{n-j},\widehat{W}^{0,x}_{t_{n-j}}),\partial_x u (t_{n-j},\widehat{W}^{0,x}_{t_{n-j}})).$$  

(ii) For $i \le n-r$, the  computation of $(Y_i,Z_i)$ involves $r$ steps  and is given by
\begin{align*} 
(PC)_r \left \{
\begin{array}{rcl}
\whZi{}  &=&\whEFp{\ti{}}{ 
		\whHi{,1} \whYi{+1} + h_i \sum_{j = 1}^r \whHi{,j}  \sum_{j=1}^r \beta_{j}  f(\whYi{+j},\whZi{+j})  }
\vspace{2mm}
\\
{}^p\whYi{} &=& \whEFp{\ti{}}{ \whYi{+1} + \hi \sum_{j=1}^r \beta_{i,j,r} f(\whYi{+j},\whZi{+j}) }
\vspace{2mm}
\\
\whYi{} &=& \whEFp{\ti{}}{ \whYi{+1} + \hi \sum_{j=1}^r b_{i,j,r} f(\whYi{+j},\whZi{+j}) } + \hi b_{i,0,r} f({}^p\whYi{},\whZi{}) 
\end{array}
\right. 
\end{align*}
where the $b$-coefficients are given by \eqref{eq de b AM}, the $\beta$-coefficients by \eqref{eq de tilde L} and the $\widehat{H}$-coefficients by \eqref{de whH}.
\end{Definition}

\begin{Proposition} \label{pr conv results}
(i) In Definition \ref{de lms numerical scheme}, if we assume that the method  given by the coefficients $a$, $b$, $\alpha$, $\beta$ is of order $m$, according to Definition \ref{de order}, and that the multinomial approximation of the Brownian motion is of order $K = 2m+1$  then we have
\begin{align*}
Y_0 - \wh{Y}_0 = O(h^m)\;,
\end{align*}
provided that the coefficient $f$ and the value function $u$ are smooth enough.

\vspace{1mm}
(ii) In Definition \ref{de PC numerical scheme} for $(PC)_r$ method, if we assume that  the multinomial approximation of the Brownian motion is of order $K = 2r+3$  then we have
\begin{align*}
Y_0 - \wh{Y}_0 = O(h^{r+1})\;,
\end{align*}
provided that the coefficient $f$ and the value function $u$ are smooth enough.
\end{Proposition}
The proof of this proposition is postponed to the end of this section.

\vspace{2mm}
We can now turn to a concrete example which illustrates the above order of convergence.
\subsection{Application}
As in \cite{zhazha10}, we consider the  process, on $[0,T]$,
\begin{align*}
(X_t,Y_t,Z_t) = \Big(W_t,\frac1{1 + \exp(-W_t-\frac{t}{4})},\frac{\exp(-W_t-\frac{t}4)}{(1+\exp(-W_t-\frac{t}4))^2}\Big) \;.
\end{align*}

This process is solution of the (decoupled) FBSDE
\begin{align*}
X_t &= W_t
\\
Y_t &= g_T(W_T) + \int_t^Tf(Y_s,Z_s) \ud s - \int_t^T Z_s \ud W_s
\end{align*}
where the driver $f$ is given by
\begin{align}
f(y,z) = -z(\frac34-y)\;, \label{de num fyz}
\quad\text{ and }\quad
g_T(x) = \frac1{1 + \exp(-x-\frac{T}{4})}
\end{align}

\newpage
To approximate the value of $Y_0$, we consider the following methods:
\begin{enumerate}
\item  Implicit Euler approximation, coupled with an order 3 Brownian approximation. 
\item Crank-Nicholson approximation, coupled with an order 5 Brownian approximation. 

\item Explicit two step Adams method, coupled with an order 5 Brownian approximation.
\item Implicit two step Adams method, coupled with an order 7 Brownian approximation. 

\item  Heun method which is a Predictor-Corrector method, coupled with an order 5 Brownian approximation.

\end{enumerate}

The log-log graph in Fig$.$ 1 below shows the rates of convergence of the method which are in accordance with the theoretical ones. Adams methods produce empirical rate slightly below the expected ones. But the highest is the theoritical convergence order, the smallest is the error in practice.

\begin{figure}[h!] \label{fig1}
\includegraphics[width = \textwidth]{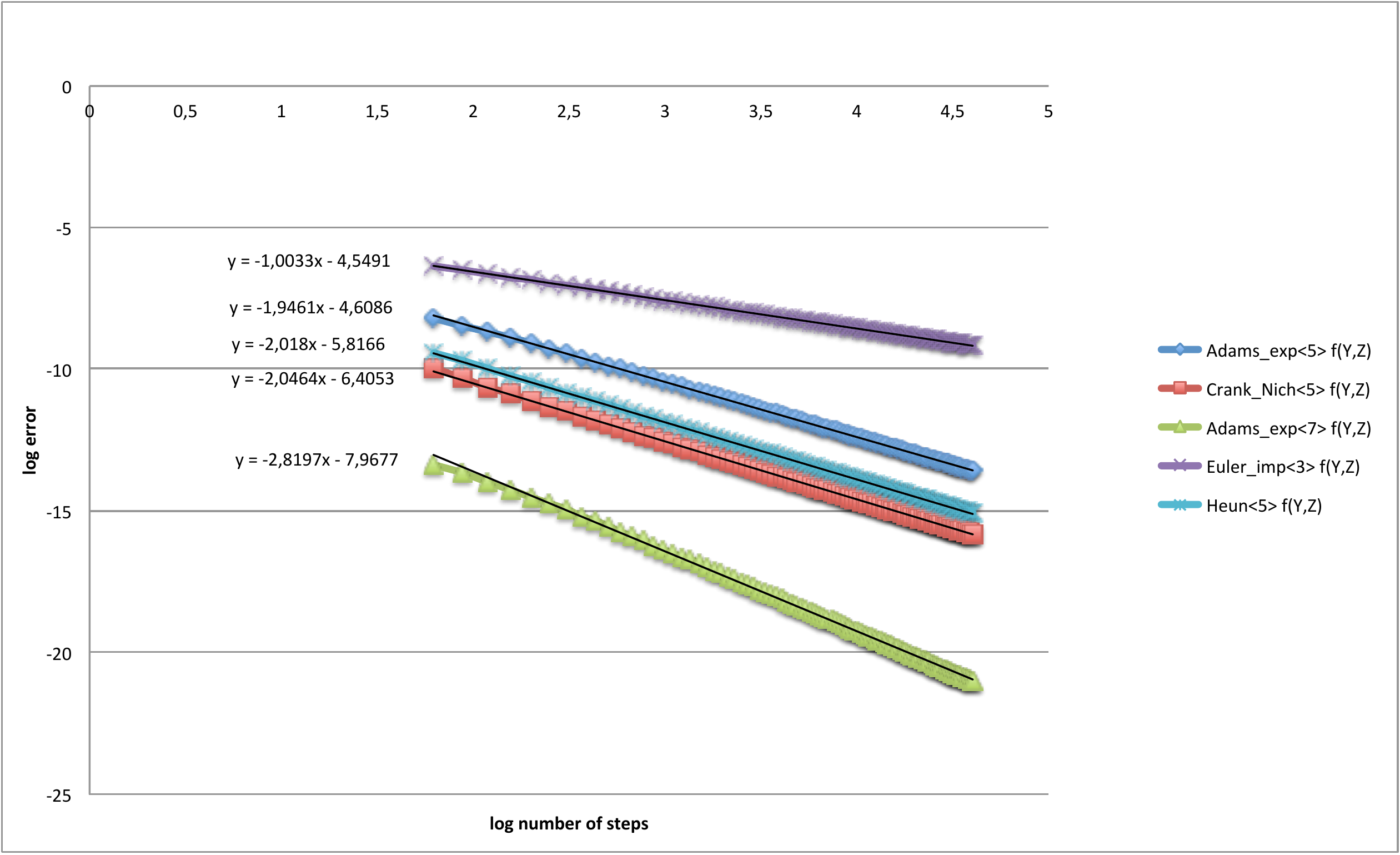}
\caption{Illustration of the convergence rate}
\end{figure}

The  graph in Fig$.$ 2 below  shows the impact of the space discretization on the global order of the method. The empirical convergence rates are in accordance with the theoretical ones.

\vspace{1mm}
\begin{figure}[h!] \label{fig2}
\includegraphics[width = \textwidth]{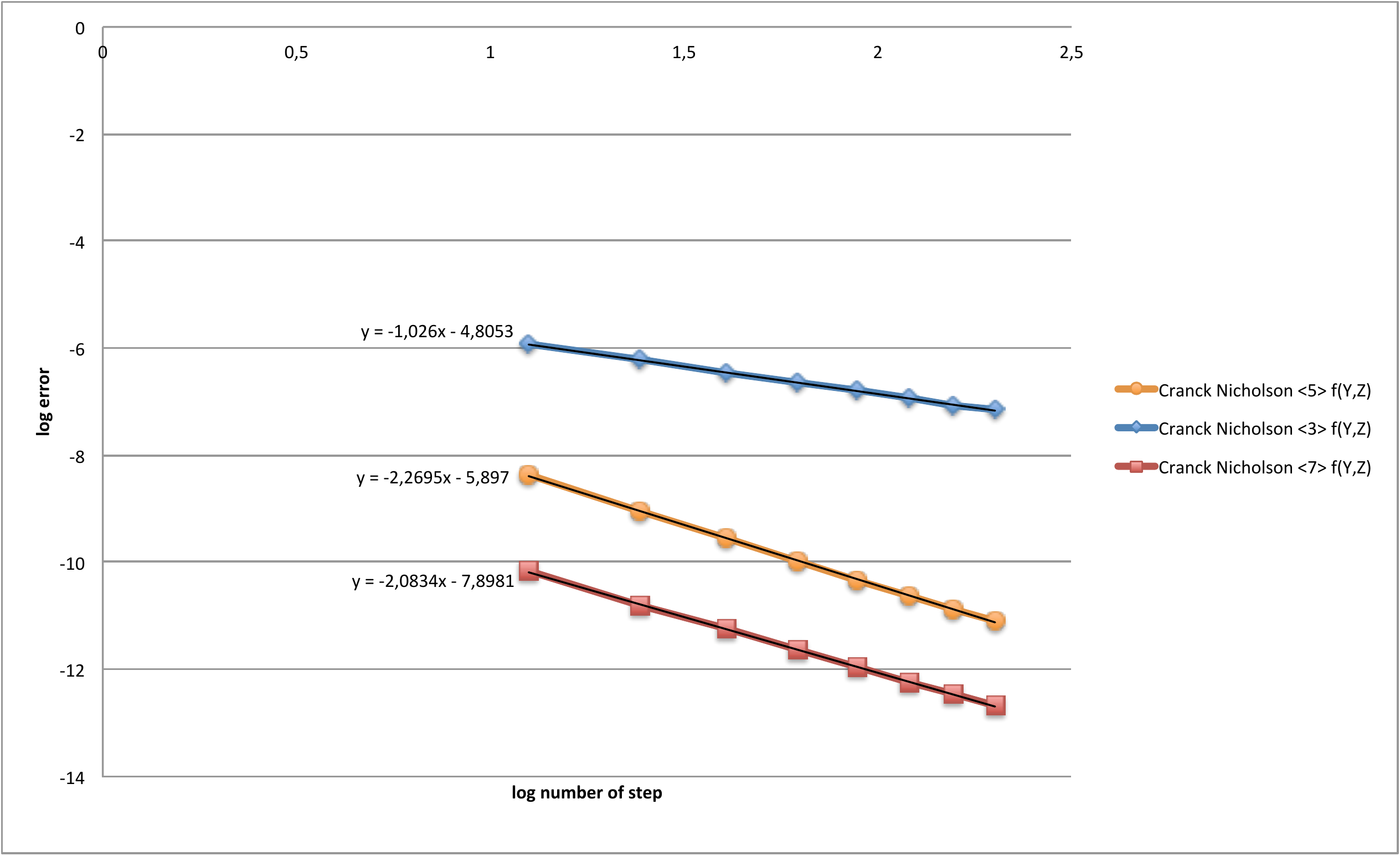}
\caption{Impact of space discretization}
\end{figure}

\subsection{Proof of Proposition \ref{pr conv results}}
We only provide the proof of (i), the proof of (ii) follows from the same arguments and using the proof of Proposition \ref{pr pc stab} and Proposition \ref{pr pc truncation error}.

\vspace{1mm}
1. \emph{Notations}

\vspace{1mm}
We first need to consider 'functional' version of the schemes above.
Let us introduce the following operator, related to  the theoretical schemes given in Definition \ref{de multi-step scheme}. 

$\Rzi{j}:(C^1_b)^{2}\rightarrow C^1_b$
\begin{align*}
\Rzi{j}[\phiy,\phiz](x) = \esp{H^{\mathbf{1},x}_{\ti{},jh} \Big( \alpha_{j}  \phiy(\Wxi{+j}) + h  \beta_{j} f(\phiy(\Wxi{+j}),\phiz(\Wxi{+j})) \Big)}
\end{align*}
$\Ryi{j}:(C^1_b)^{2}\rightarrow C^1_b$
\begin{align*}
\Ryi{j}[\phiy,\phiz](x) = \esp{a_{j}  \phiy(\Wxi{+j}) + h  b_{j} f(\phiy(\Wxi{+j}),\phiz(\Wxi{+j})) }
\end{align*}

%
%
%
%
\vspace{2mm}
Similarly, let us define - for the fully discrete scheme - the operators

$\whRzi{j}:(C^1_b)^{2}\rightarrow C^1_b$
\begin{align*}
\whRzi{j}[\phiy,\phiz](x) = \whesp{ \widehat{H}_{i,j} \Big( \alpha_{j}  \phiy(\whWxi{j}) + h  \beta_{j} f(\phiy(\whWxi{j}),\phiz(\whWxi{j})) \Big)}
\end{align*}
$\whRyi{j}:(C^1_b)^{2}\rightarrow C^1_b$
\begin{align*}
\whRyi{j}[\phiy,\phiz](x) = \whesp{a_{j}  \phiy(\whWxi{+j}) + h  b_{j} f(\phiy(\whWxi{+j}),\phiz(\whWxi{+j})) }
\end{align*}

The functional version of the  schemes given in Definition \ref{de lms numerical scheme} reads then, for $i \le n-r$,
\begin{align*} 
\left \{
\begin{array}{rcl}
\whyi{}(x) &=& \sum_{j=1}^r \whRyi{j}[\whyi{+j},\whzi{+j}](x)
\vspace{2mm}
\\
\whzi{}(x) &=& \sum_{j=1}^r \whRzi{j}[\whyi{+j},\whzi{+j}](x)
\end{array}
\right.
\end{align*}

given $r$ initial data $(\widehat{y}_{n-j},\widehat{z}_{n-j})=(u(t_{n-j},\cdot), \partial_x u(t_{n-j},\cdot))$, $0\le j \le r-1$. 

\vspace{2mm}
Due to the markov property of the discrete process $(\widehat{W}^{0,x}_{t})_{t\in \pi}$, it is easily checked that
\begin{align*}
\whYi{} = \whyi{}(\Wpi{}) \text{ and } \whZi{} = \whzi{}(\Wpi{}) \;.
\end{align*}
Finally, we define
\begin{align*}
\wtYi{} = u(\ti{},\Wpi{}) \text{ and } \wtZi{} = \partial_x u(\ti{},\Wpi{}) \;.
\end{align*}
Observe that $\widetilde{Y}_0 = u(0,x)$ and that, for $0\le j \le r-1$, $(\widehat{Y}_{n-j},\widehat{Z}_{n-j}) = (\widetilde{Y}_{n-j},\widetilde{Z}_{n-j})$.

\vspace{2mm}
2. \emph{Stability}

\vspace{1mm}
The key observation is here that $(\wtYi{},\wtZi{})$ can be seen as a perturbed version of the scheme given in \eqref{de numerics scheme}, namely
\begin{align*}
\left \{
\begin{array}{rcl}
\wtYi{} &=& \whEFp{\ti{}}{ \sum_{j=1}^r a_{j}\wtYi{+j}+ h \sum_{j=0}^r b_{j} f(\wtYi{+j},\wtZi{+j}) }  + {}^{t}\errYi{} + {}^{s}\errYi{}
\vspace{2mm}
\\
\wtZi{}  &=&\whEFp{\ti{}}{ 
		\sum_{j = 1}^r H \Big( \alpha_{j} \wtYi{+j} + h \sum_{j=1}^r \beta_{j}  f(\wtYi{+j},\wtZi{+j})  \Big) }+ {}^{t}\errZi{} + {}^{s}\errZi{}
\end{array}
\right.
\end{align*}
where the local error due to the time-discretization is
\begin{align}
\left \{
\begin{array}{rcl}
{}^t\errYi{} &=&\esp{Y_{\ti{}} - \check{Y}_{\ti{}} | X_{\ti{}} = \Wpi{}}
\vspace{2mm}
\\
{}^t\errZi{} &=&  \esp{Z_{\ti{}} - \check{Z}_{\ti{}} | X_{\ti{}} = \Wpi{}}
\end{array}
\right.
\end{align}
recalling \eqref{eq def hY}-\eqref{eq def hZ}
and the local error due to the 'space-discretization' is
\begin{align}
\left \{
\begin{array}{rcl}
{}^s\errYi{} &=& \sum_{j=1}^r (\Ryi{j} - \whRyi{j})[u(\ti{+j},\cdot),\partial_x u(\ti{+j},\cdot)](\Wpi{})
\vspace{1mm}
\\
{}^s\errZi{} &=&  \sum_{j=1}^r (\Rzi{j} - \whRzi{j})[u(\ti{+j},\cdot),\partial_x u(\ti{+j},\cdot)](\Wpi{})
\end{array}
\right. .
\end{align}

Now, we can apply Proposition \ref{pr L2 stab}, recalling Remark \ref{re pr L2 stab}, to obtain in particular that
\begin{align} \label{eq stab numerics}
|\widetilde{Y}_0 - \widehat{Y}_0|^2 \le  |\pi|\sum_{i=0}^{n-r} \whesp{\frac{1}{\hi^2} |{}^t\errYi{} +{}^s\errYi{}|^2 +  |{}^t\errZi{} +{}^s\errZi{}|^2}
\end{align}

\vspace{2mm}
3. \emph{Study of the local error}

\vspace{1mm}
We now turn to the study of the local errors $(\errYi{},\errZi{})_{0 \le i \le n-r}$.
Assuming that the function are smooth enough we compute the following expansion
\begin{align*}
(\Ryi{j} - \whRyi{j})[u(\ti{+j},\cdot),u^{(1)}(\ti{+j},\cdot)](x) = 
\sum_{k=0}^K \frac1{k!} \chi^{(k)}_{i,j}(x)\esp{(\Woi{+j})^k} &- \sum_{k=0}^K \frac1{k!} \chi^{(k)}_{i,j}(x) \whesp{(\whWoi{+j})^k}  
\\&+ O(|\pi|^{\frac{K+1}{2}})
\end{align*}
where $\chi_{i,j}$ are functions depending on $f$, $u$ and the coefficients of the methods.

Using the matching moment property of $\widehat{W}^{0,\ti{}}$, we easily obtain that 
\begin{align*}
\whesp{|(\Ryi{j} - \whRyi{j})[u(\ti{+j},\cdot),u^{(1)}(\ti{+j},\cdot)](\Wpi{})|^2} \le C|\pi|^{K+1}
\end{align*}

For the $Z$ part, we have
\begin{align*}
(\Rzi{j} - \whRzi{j})[u(\ti{+j},\cdot),u^{(1)}(\ti{+j},\cdot)](x) = \sum_{k=0}^{K-1} \frac1{k!} \chi^{(k)}_{i,j}(x)\esp{H_{i,j}(\Woi{+j})^k} &- \sum_{k=0}^{K-1} \frac1{k!} \chi^{(k)}_{i,j}(x) \whesp{(\widehat{H}_{i,j}\whWoi{+j})^k}  \\
& +O(|\pi|^\frac{K-1}{2})
\end{align*}
Using the matching moment property of $\widehat{W}^{0,\ti{}}$, we easily obtain that 
\begin{align*}
\whesp{|(\Rzi{j} - \whRzi{j})[u(\ti{+j},\cdot),u^{(1)}(\ti{+j},\cdot)](\Wpi{})|^2} \le C |\pi|^{K-1}
\end{align*}

Combining the above estimates with \eqref{eq stab numerics} and the fact that the discrete-time error is of order $m$, leads to
\begin{align*}
|\widetilde{Y}_0 - \widehat{Y}_0| \le C|\pi|^{m} + C|\pi|^{\frac{K-1}{2}}.
\end{align*}
which concludes the proof since $K=2m+1$. \eproof

\renewcommand{\ti}{t_{i}}

\end{document}